\documentclass[review]{elsarticle}
\usepackage{url}
\usepackage{amssymb}
\usepackage{amsthm}
\theoremstyle{plain}

\theoremstyle{remark}
\newtheorem{rem}{\protect\remarkname}
\theoremstyle{plain}
\newtheorem{obs}{\protect\observationname}
\theoremstyle{definition}
\newtheorem{defn}{\protect\definitionname}
\theoremstyle{prop}

\theoremstyle{definition}

\theoremstyle{plain}

\theoremstyle{plain}

\providecommand{\definitionname}{Definition}
\providecommand{\propositionname}{Proposition}
\providecommand{\examplename}{Example}
\providecommand{\lemmaname}{Lemma}
\providecommand{\observationname}{Observation}
\providecommand{\remarkname}{Remark}
\providecommand{\theoremname}{Theorem}
\providecommand{\corollaryname}{Corollary}
\usepackage{graphicx}
\usepackage{subfig}
\usepackage[cmex10]{amsmath}
\usepackage{color}
\usepackage{algorithm,algorithmicx}
\usepackage{caption}
\usepackage{appendix}
\makeatletter
\@addtoreset{algorithm}{chapter}
\makeatother
\newcommand\numberthis{\addtocounter{equation}{1}\tag{\theequation}}

\def \X {\mathbf{X}}

\def \C {\mathbf{C}}

\def \x {\mathbf{x}}

\def \s {\mathbf{s}}

\def \v {\mathbf{v}}

\begin{document}
\begin{frontmatter}{}
\title{Simultaneous diagonalisation of the covariance and complementary covariance matrices in quaternion widely linear signal processing}
\author[address1]{Min Xiang\corref{correspondingauthor}}
\cortext[correspondingauthor]{Corresponding author}
\ead{m.xiang13@ic.ac.uk}
\author[address2]{Shirin Enshaeifar}
\author[address1]{Alexander E. Stott}
\author[address2]{Clive Cheong Took}
\author[address3]{Yili Xia}
\author[address1]{Sithan Kanna}
\author[address1]{Danilo P. Mandic}
\address[address1]{Department
of Electrical and Electronic Engineering, Imperial College London, London SW7 2BT, UK}
\address[address2]{Department of Computer Science, the University of Surrey, Surrey GU2 7XH, UK}
\address[address3]{School of Information Science and Engineering, Southeast University, Nanjing 210018, China}


\begin{abstract}
Recent developments in quaternion-valued widely linear processing have established that the exploitation of complete second-order statistics requires consideration of both the standard covariance and the three complementary covariance matrices. Although such matrices have a tremendous amount of structure and their decomposition is a powerful tool in a variety of applications, the non-commutative nature of the quaternion product has been prohibitive to the development of quaternion uncorrelating transforms. To this end, we introduce novel techniques for a simultaneous decomposition of the covariance and complementary covariance matrices in the quaternion domain, whereby the quaternion version of the Takagi factorisation is explored to diagonalise symmetric quaternion-valued matrices. This gives new insights into the quaternion uncorrelating transform (QUT) and forms a basis for the proposed quaternion approximate uncorrelating transform (QAUT) which simultaneously diagonalises all four covariance matrices associated with improper quaternion signals. The effectiveness of the proposed uncorrelating transforms is validated by simulations on both synthetic and real-world quaternion-valued signals.
\end{abstract}
\begin{keyword}
Quaternion matrix diagonalisation, complementary covariance, pseudo-covariance, widely linear processing, quaternion non-circularity, improperness.
\end{keyword}
\end{frontmatter}{}

\section{Introduction}
Advances in sensing technology have enabled widespread recording from 3-D and 4-D data sources, such as measurements from seismometers \cite{Krieger2015}, ultrasonic anemometers \cite{CheongTook2011c}, and inertial body sensors \cite{Fourati2014}. It is convenient to express such measurements as vectors in the $\mathbb{R}^3$ and $\mathbb{R}^4$ fields of reals, however, the real vector algebra is not a division algebra and is therefore inadequate when modelling orientation and rotation \cite{Kuipers1999}. Quaternions have advantages in representing 3-D and 4-D data, owing to their division algebra, a generic extension of the real and complex algebras. Quaternions also naturally account for mutual information between multiple data channels, provide a compact representation, and have proven to offer a physically meaningful interpretation for real-world applications in the fields of navigation, communication, and image processing \cite{Tao2014,Chen2014}. A recent resurgence in research on quaternion signal processing spans the areas of filtering \cite{CheongTook2009}, independent component analysis (ICA) \cite{Via2011,Javidi2011}, neural networks \cite{Xia2015,Minemoto2017}, and Fourier transforms \cite{Ell2014}.

Diagonalisable matrices play a fundamental role in engineering applications, whereby the related uncorrelating transforms greatly simplify the analysis of complex problems, both in terms of enabling a solution and providing a methodological framework (e.g. performance bounds). Computational advantages associated with matrix diagonalisation include a reduction in the size of the parameter space from $\mathbf{N}^2$ to $\mathbf{N}$, while the statistical advantages include the possibility to match source properties (e.g. orthogonality) in blind source separation. In addition, joint diagonalisation of matrices in some blind source separation applications has a close link with the canonical decomposition problem for tensors \cite{DeLathauwer2006,Cichocki2015}. The real linear algebra is a mature area, however, covariance matrices in widely linear signal processing \cite{de2002} in the complex and quaternion division algebras and their structures \cite{Eriksson2006} have recently received significant attention \cite{schreier2010statistical,adali2010adaptive,yeredor2012,Via2010,CheongTook2011b}. 

In the context of complex-valued signal processing, for a complex random vector $\mathbf{z}$, both its covariance, $\mathbf{C}=E\{\mathbf{zz}^H\}$, and pseudo-covariance, $\mathbf{P}=E\{\mathbf{zz}^T\}$, matrices are necessary to capture complete second-order statistical information \cite{Mandic2009}. To extract the statistical information embedded in the two covariance matrices, simultaneous diagonalisation of $\mathbf{C}$ and $\mathbf{P}$ is a prerequisite, while their respective individual diagonalisations can be performed via the eigendecomposition and the Takagi factorisation. Such decorrelations exist in their general forms in mathematics \cite{horn1990matrix}. De Lathauwer and De Moor \cite{de2002} and Eriksson and Koivunen \cite{Eriksson2006} have given them a practical value through an important engineering contribution called the strong uncorrelating transform (SUT), together with its extension, the generalised uncorrelating transform \cite{Ollila2009}. Our own recent contributions include a computationally efficient approximate uncorrelating transform (AUT) \cite{CheongTook2012}, while to preserve the integrity of the original bivariate sources in the decorrelation process, we have proposed the correlation preserving transform \cite{CheongTook2015}. These transforms have been used in performance analysis of complex-valued adaptive filters \cite{douglas2010performance,mandic2015mean,xia2017complementary,xia2017full}, blind separation of non-circular complex-valued sources \cite{shen2010complex,Ramírez2012}, complex-valued subspace tracking \cite{Douglas2006}, and separation of signal and noise components in subbands of harmonic signals \cite{Okopal2015}.

Real- and complex-valued signal processing have benefited greatly from making sophisticated use of real- and complex-valued matrix algebras \cite{bojanczyk,hjorungnes2011complex}. However, quaternion-valued signal processing, which is advantageous for 3-D and 4-D data, has been hindered by the underdevelopment of quaternion-valued matrix algebra \cite{Zhang1997}, compared to the well-established real- and complex-valued matrix algebras \cite{horn1990matrix,hjorungnes2011complex}. This is due to the non-commutativity nature of quaternion multiplication and the enhanced degrees of freedom provided by quaternions. For example, a quaternion square matrix has right and left eigenvalues; the right eigenvalues have been well studied, while the left eigenvalues are less known and are not computationally well-posed \cite{Zhang1997,zou2012location}. The recent theory of quaternion matrix derivatives provides a systematic framework for the calculation of derivatives of quaternion matrix functions with respect to quaternion matrix variables \cite{Xu2015a}. Advances in structural quaternion matrix decompositions include tools to diagonalise quaternion matrices \cite{Sangwine2006, Bihan2007}, while in the context of quaternion widely linear processing, it was shown that for a quaternion random vector $\mathbf{x}$, the diagonalisation of the Hermitian covariance matrix $\mathbf{C}_{\mathbf{x}}=E\{\mathbf{\mathbf{xx}}^{H}\}$ can be performed straightforwardly using the quaternion eigendecomposition, whereas the diagonalisation of the complementary covariance matrices ($\mathbf{C}_{\mathbf{x}^\imath}=E\{\mathbf{\mathbf{xx}}^{\imath H}\},~\mathbf{C}_{\mathbf{x}^\jmath}=E\{\mathbf{\mathbf{xx}}^{\jmath H}\},~\mathbf{C}_{\mathbf{x}^\kappa}=E\{\mathbf{\mathbf{xx}}^{\kappa H}\}$) can be computed via the quaternion singular value decomposition (SVD) \cite{CheongTook2011a}. Furthermore, the quaternion uncorrelating transform (QUT), a quaternion analogue of the SUT for complex data, was proposed in \cite{ShirinICASSP} to jointly diagonalise the covariance matrix and one of the three complementary covariance matrices. However, the simultaneous diagonalisation of all the four covariance matrices is still an open problem. It is important to notice that the relationship between the quaternion covariance matrices in the widely linear model is governed by \cite{CheongTook2011b} 
\begin{equation} \mathbf{P}_\mathbf{x}=\frac{1}{2}(\mathbf{C}_{\mathbf{x}^\imath}+\mathbf{C}_{\mathbf{x}^\jmath}+\mathbf{C}_{\mathbf{x}^\kappa}-\mathbf{C}_{\mathbf{x}})\label{eq: covariance relationship}
\end{equation}
which demonstrates that the diagonalisation of the pseudo-covariance matrix, $\mathbf{P}_\mathbf{x}=E\{\mathbf{\mathbf{xx}}^{T}\}$, requires us to address the following issues: 
\begin{itemize}
\item There is no \emph{closed-form} solution\footnote{Iterative solutions were proposed in the context of quaternion ICA \cite{Via2011}.} to perform a simultaneous diagonalisation of the covariance matrix and the three complementary covariance matrices. 
\item Given the dimensionality of the problem and numerous practical applications, a simple approximate diagonalising transform which would apply to both the covariance and pseudo-covariance matrices would be beneficial. 
\end{itemize}
In this paper, we therefore set out to propose solutions to these open problems, and support the analysis with illustrative examples.

The rest of this paper is organised as follows. Section \ref{sect: background} provides an overview
of quaternion algebra and statistics. Section \ref{sect: diagonalisation} proposes novel techniques for the diagonalisation of symmetric quaternion matrices and a simultaneous diagonalisation of two quaternion covariance matrices. Section \ref{sect:QAUT} proposes an approximate simultaneous diagonalisation of the four quaternion covariance matrices. Simulation results are given in Section \ref{sect: simulation}, and Section \ref{sect: conclusion} concludes the paper. Throughout the paper, we use
boldface capital letters to denote matrices, $\mathbf{A}$, boldface
lowercase letters for vectors, $\mathbf{a}$, and standard letters
for scalar quantities, $a$. Superscripts $(\cdot)^{T}$, $(\cdot)^{*}$
and $(\cdot)^{H}$ denote the transpose, conjugate, and Hermitian
(i.e. transpose and conjugate), respectively, $\mathbf{I}$ the identity matrix,
and $E\left\{ \cdot\right\} $ the statistical expectation operator.
\section{Background}\label{sect: background}
\subsection{Quaternion algebra}
The quaternion domain $\mathbb{H}$ is a 4-D vector space over the
field of reals, spanned by the basis $\left\{ 1,\imath,\jmath,\kappa\right\}$, where $\imath=(1,0,0)$, $\jmath=(0,1,0)$ and $\kappa=(0,0,1)$ are three unit vectors in $\mathbb{R}^3$.
A quaternion vector\footnote{Throughout this paper, the quaternion vector is assumed to be zero-mean.} This does not affect the generality of our results. $\mathbf{x}\in \mathbb{H}^{L\times1}$ consists of a real (scalar) part $\mathfrak{R}[\cdot]$ and an imaginary (vector) part $\mathfrak{I}[\cdot]$ which comprises $\imath-$, $\jmath-$, and $\kappa-$ imaginary components, so that
\begin{equation}
\begin{split}
\mathbf{x}=&\mathfrak{R}[\mathbf{x}]+\mathfrak{I}[\mathbf{x}]\\
    =&\mathfrak{R}[\mathbf{x}]+\mathfrak{I}_\imath[\mathbf{x}]\imath+\mathfrak{I}_\jmath[\mathbf{x}]\jmath+\mathfrak{I}_\kappa[\mathbf{x}]\kappa \\  =&\mathbf{x}_a+\mathbf{x}_b\imath+\mathbf{x}_c\jmath+\mathbf{x}_d\kappa
\end{split} \label{eq:quatvector}
\end{equation}
\noindent where $\mathfrak{R}[\mathbf{x}]=\mathbf{x}_a, \mathfrak{I}_\imath[\mathbf{x}]=\mathbf{x}_b, \mathfrak{I}_\jmath[\mathbf{x}]=\mathbf{x}_c, \mathfrak{I}_\kappa[\mathbf{x}]=\mathbf{x}_d$ are $L\times1$ real vectors, and $\imath$, $\jmath$, $\kappa$ are orthogonal imaginary units with properties
\begin{equation}
\begin{split}
\imath\jmath&=-\jmath\imath=\kappa  \qquad \jmath\kappa=-\kappa\jmath=\imath \qquad \kappa\imath=-\imath\kappa=\jmath \\
\imath^{2}&=\jmath^{2}=\kappa^{2}=-1
\end{split} \label{eq:units}
\end{equation}
A quaternion variable $x$ is called a pure quaternion if it satisfies $\mathfrak{R}[x]=0$. The modulus of a quaternion variable $x\in\mathbb{H}$ is defined as \[|x|=\sqrt{x_{a}^{2}+x_{b}^{2}+x_{c}^{2}+x_{d}^{2}}\]
A quaternion variable $x$ is called a unit quaternion if it satisfies $|x|=1$. 
The product of two quaternions $x,y\in\mathbb{H}$ are defined as 
\[
xy=\mathfrak{R}[x] \mathfrak{R}[y]-\mathfrak{I}[x]\cdot\mathfrak{I}[y]+\mathfrak{R}[x]\mathfrak{I}[y]+\mathfrak{R}[y]\mathfrak{I}[x]+\mathfrak{I}[x]\times\mathfrak{I}[y]
\] 
where ``$\cdot$'' denotes the scalar product and ``$\times$'' the vector product. Notice that the presence of the vector product in the above expression causes the non-commutativity of the quaternion product, that is, $xy\neq yx$.

Another important notion in the quaternion domain is the so-called ``quaternion involution'' \cite{Ell2007}, which defines a self-inverse mapping analogous to the complex conjugate. The general involution of the quaternion vector $\mathbf{x}$ is defined as $\mathbf{x}^{\alpha}=-\alpha \mathbf{x}\alpha$, and represents the rotation of the vector part of $\mathbf{x}$ by $\pi$ about a unit pure quaternion $\alpha$. The special cases of involutions about the $\imath$, $\jmath$ and $\kappa$ imaginary axes are given by~\cite{Ell2007}
\begin{equation}
\begin{split}
\mathbf{x}^\imath &= -\imath\mathbf{x}\imath = \mathbf{x}_a + \mathbf{x}_b \imath- \mathbf{x}_c \jmath- \mathbf{x}_d \kappa\\
\mathbf{x}^\jmath &= -\jmath\mathbf{x}\jmath = \mathbf{x}_a - \mathbf{x}_b \imath+ \mathbf{x}_c \jmath- \mathbf{x}_d \kappa\\
\mathbf{x}^\kappa &= -\kappa\mathbf{x}\kappa = \mathbf{x}_a - \mathbf{x}_b \imath-\mathbf{x}_c \jmath + \mathbf{x}_d\kappa
\end{split}\label{eq: invol}
\end{equation}
This set of involutions, together with the original quaternion, $\mathbf{x}$, forms the most frequently used basis for augmented quaternion statistics \cite{CheongTook2011b,Via2010}. While an involution represents a rotation along a single imaginary axis, the quaternion conjugate operator $(\cdot)^*$ rotates the quaternion along all three imaginary axes, and is given by
\begin{equation}
\mathbf{x}^{*}=\mathfrak{R}[\mathbf{x}]-\mathfrak{I}[\mathbf{x}]= \mathbf{x}_a - \mathbf{x}_b \imath- \mathbf{x}_c \jmath- \mathbf{x}_d\kappa
\label{eq: conjugate}
\end{equation}
A useful algebraic operation which conjugates only one imaginary component is given by
\begin{equation}
\begin{split}
\mathbf{x}^{\imath*}&=(\mathbf{x}^{\imath})^*=(\mathbf{x}^{*})^\imath = \mathbf{x}_a - \mathbf{x}_b\imath + \mathbf{x}_c\jmath + \mathbf{x}_d \kappa\\
\mathbf{x}^{\jmath*}&= (\mathbf{x}^{\jmath})^*=(\mathbf{x}^{*})^\jmath = \mathbf{x}_a + \mathbf{x}_b\imath - \mathbf{x}_c\jmath + \mathbf{x}_d\kappa \\
\mathbf{x}^{\kappa*}&= (\mathbf{x}^{\kappa})^*=(\mathbf{x}^{*})^\kappa = \mathbf{x}_a + \mathbf{x}_b\imath + \mathbf{x}_c\jmath - \mathbf{x}_d\kappa
\end{split}\label{eq: invol_conj}
\end{equation}
Since the only difference between $\mathbf{x}$ and $\mathbf{x}^{\alpha*}$, for $\alpha\in\{\imath,\jmath,\kappa\}$, is the sign of the $\alpha$-imaginary component, a quaternion vector $\mathbf{x}$ is called $(\cdot)^{\alpha*}$ invariant if $\mathfrak{I}_\alpha[\mathbf{x}]=0$. For example, if $\mathbf{x}$ has a vanishing $\imath$-imaginary part, then $\mathbf{x}^{\imath *}=\mathbf{x}$.
For rigour, the structure and general properties of quaternion matrices are given in Appendix A.
\subsection{Quaternion covariance matrices: Structural insights}
The recent success of complex augmented statistics is largely due to the simplicity and physical meaningfulness of the statistical descriptors in the form of the standard covariance, $\mathbf{C}_\mathbf{x}=E\{\mathbf{xx}^H\}$, and the pseudo-covariance, $\mathbf{P}_\mathbf{x}=E\{\mathbf{xx}^T\}$. In particular, the pseudo-covariance enables us to account for the improperness (channel power imbalance or correlation) of complex variables, while the symmetry of $\mathbf{P}_\mathbf{x}$ implies that its decomposition can be computed using the Takagi factorisation as $\mathbf{P}_\mathbf{x}=\mathbf{Q}\boldsymbol{\Sigma}\mathbf{Q}^T$, where $\mathbf{Q}$ is a complex unitary matrix and $\boldsymbol{\Sigma}$ a real non-negative diagonal matrix \cite{horn1990matrix}. However, the pseudo-covariance of a quaternion random vector $\mathbf{x}=[x_1,x_2,\ldots,x_L]^T$ does not exhibit symmetry, owing to the non-commutativity of the quaternion product, whereby $x_mx_n \neq x_nx_m\;(m\neq n)$, yield an asymmetric matrix
\begin{equation}
\begin{split}
\!\mathbf{P}_{\mathbf{x}}\!&\!=\!E\{\mathbf{xx}^{T}\}\!\!=\!\!\begin{bmatrix}
{E\{x_1x_1\}} & {E\{x_1x_2\}} & \ldots & {E\{x_1x_L\}}\\
{E\{x_2x_1\}} & {E\{x_2x_2\}} & \ldots & {E\{x_2x_L\}}\\
\vdots & \vdots & \ddots & \vdots  \\
{E\{x_Lx_1\}} & {E\{x_Lx_2\}} & \ldots & {E\{x_Lx_L\}}
\end{bmatrix}
\end{split}
\end{equation}

The quaternion involution basis in \eqref{eq: invol} is at the core of the recently proposed widely linear processing \cite{Via2010,CheongTook2010}, which provides theoretical and practical performance gains over traditional strictly linear processing \cite{Xia2015a}. For example, the widely linear minimum mean square error (MMSE) estimation of the quaternion signal $\mathbf{y}$ in terms of the observation $\mathbf{x}$ is performed as
\begin{equation}\label{WL}
\mathbf{\hat{y}}=E\{\mathbf{y}|\mathbf{x},\mathbf{x}^\imath,\mathbf{x}^\jmath,\mathbf{x}^\kappa\}
\end{equation}
which is analogous to the complex widely linear MMSE estimation \cite{Mandic2009}, given by
\[
\mathbf{\hat{y}}=E\{\mathbf{y}|\mathbf{x},\mathbf{x}^*\}
\]
Then, the involution basis for the quaternion widely linear processing provides a useful description of quaternion second-order statistics, represented by the standard covariance matrix and the $\imath-$, $\jmath-$, and $\kappa-$ complementary covariance matrices, which are given by \cite{CheongTook2011b}
\begin{equation}
{\small
\begin{split}
\mathbf{C}_{\mathbf{x}^\alpha}\!\!=\!\!E\{\mathbf{x}\mathbf{x}^{\alpha H}\}\!\!=\!\!\begin{bmatrix}
E\{x_1x_1^{\alpha *}\} & E\{x_1x_2^{\alpha *}\} & \cdots & E\{x_1x_L^{\alpha *}\}\\
E\{x_2x_1^{\alpha *}\} & E\{x_2x_2^{\alpha *}\} & \cdots & E\{x_2x_L^{\alpha *}\} \\
\vdots & \vdots & \ddots & \vdots  \\
E\{x_Lx_1^{\alpha *}\} & E\{x_Lx_2^{\alpha *}\} & \cdots & E\{x_Lx_L^{\alpha *}\}
\end{bmatrix} 
\end{split}}
\end{equation}
where $\alpha\in\{\imath,\jmath,\kappa\}$. The $\alpha$-complementary covariance matrix is $\alpha$-Hermitian, that is, $\mathbf{C}_{\mathbf{x}^\alpha}=(\mathbf{C}_{\mathbf{x}^\alpha})^{\alpha H}$, which stems from the fact that its diagonal entries are $(\cdot)^{\alpha *}$ invariant, whereas its off-diagonal entries are governed by the relationship $\mathbf{C}_{\mathbf{x}^\alpha}[m,n]=(\mathbf{C}_{\mathbf{x}^\alpha}[n,m])^{\alpha *}$, where $\mathbf{C}_{\mathbf{x}^\alpha}[m,n]$ denotes the entry of $\mathbf{C}_{\mathbf{x}^\alpha}$ in row $m$, column $n$.
The relationship between the pseudo-covariance and the three complementary covariances is given in (\ref{eq: covariance relationship}).
\begin{rem}
The knowledge of both the standard covariance matrix and the three complementary covariance matrices is necessary to ensure the exploitation of complete second-order statistical information in the quaternion domain.
\end{rem}
\subsection{Quaternion properness}
The notion of non-circularity (improperness) is unique to division algebras. While non-circularity refers to probability distributions which are not rotation-invariant, a proper complex random vector $\mathbf{z}=\mathbf{z}_r+\mathbf{z}_\imath\imath$ ($\mathbf{z}_r, \mathbf{z}_\imath \in \mathbb{R}^{L\times1}$) has a vanishing pseudo-covariance $E\{\mathbf{z}\mathbf{z}^T\}=\mathbf{0}$. In other words, its real and imaginary parts are uncorrelated and have equal variance, that is, $E\{\mathbf{z}_r\mathbf{z}_\imath^T\}=\mathbf{0}$ and $E\{\mathbf{z}_r\mathbf{z}_r^T\}=E\{\mathbf{z}_\imath\mathbf{z}_\imath^T\}$. Similarly, improperness in the quaternion domain is characterised by the degree of correlation and/or power imbalance between imaginary components relative to the real component. The additional degrees of freedom in the quaternion domain allow for two types of properness: $\mathbb{H}$-properness and $\mathbb{C}^\alpha$-properness \cite{Amblard2004}.
\begin{defn}[$\mathbb{H}$-properness]
A quaternion random vector $\mathbf{x}$ is $\mathbb{H}$-proper if it is uncorrelated with its vector involutions, $\mathbf{x}^\imath$, $\mathbf{x}^\jmath$, and $\mathbf{x}^\kappa$, so that
\begin{equation}
\begin{split}
\mathbf{C}_{\mathbf{x}^\imath}&=E\{\mathbf{x}\mathbf{x}^{\imath H}\}=\mathbf{0} \quad \quad
\mathbf{C}_{\mathbf{x}^\jmath}=E\{\mathbf{x}\mathbf{x}^{\jmath H}\}=\mathbf{0} \\
\mathbf{C}_{\mathbf{x}^\kappa}&=E\{\mathbf{x}\mathbf{x}^{\kappa H}\}=\mathbf{0}
\end{split}
\end{equation}
\end{defn}
\begin{defn}[$\mathbb{C}^\alpha$-properness]
A quaternion random vector $\mathbf{x}$ is $\mathbb{C}^\alpha$-improper\footnote{It is called $\mathbb{C}^\alpha$-proper in some literature, but we consider the term $\mathbb{C}^\alpha$-improper is more intuitive.} with respect to $\alpha=\imath,\jmath$ or $\kappa$ if it is correlated only with the involution $\mathbf{x}^\alpha$, so that all the complementary covariances except for $\mathbf{C}_{\mathbf{x}^\alpha}$ vanish.
\end{defn}
As shown in Appendix B, a $\mathbb{C}^\alpha$-improper quaternion random vector can be generated from two proper complex random vectors.
\section{Diagonalisation of quaternion covariance matrices}
\label{sect: diagonalisation}
Prior to introducing the diagonalisation of the quaternion pseudo-covariance matrix and the simultaneous diagonalisation of two quaternion covariance matrices, the following two observations establish results essential for subsequent analyses. 
\begin{obs}
The eigendecomposition applied to the Herimitian quaternion covariance matrix, $\mathbf{C}_\mathbf{x}$, gives $\mathbf{C}_\mathbf{x}=\mathbf{Q}^{H}\mathbf{\boldsymbol{\Lambda}}_{\mathbf{x}}\mathbf{Q}$, where $\mathbf{Q}$ is a quaternion unitary matrix and $\boldsymbol{\Lambda}_{\mathbf{x}}$ a real-valued diagonal matrix for which the entries are the eigenvalues\footnote{They are both right and left eigenvalues, but are simply termed eigenvalues for conciseness in this paper.} of $\mathbf{C}_\mathbf{x}$ \cite{Rodman2014}.\label{lemma2}
\end{obs}
\begin{obs}
The $\alpha$-complementary covariance matrix, $\mathbf{\mathbf{C}}_{\mathbf{x}^\alpha}$, $\alpha\in\left\{\imath,\jmath,\kappa\right\}$, can be factorised as $\mathbf{\mathbf{C}}_{\mathbf{x}^\alpha}=\mathbf{Q}_{\alpha}^{H}\mathbf{\Lambda}_{\alpha}\mathbf{Q}^{\alpha}_{\alpha}$, where $\mathbf{Q}_{\alpha}$ is a quaternion
unitary matrix and $\mathbf{\mathbf{\boldsymbol{\Lambda}}}_{\alpha}$ is a real-valued non-negative diagonal matrix for which the diagonal entries are the singular values of $\mathbf{C}_{\mathbf{x}^\alpha}$ \cite{CheongTook2011a}.\label{lemma1}
\end{obs}
Observation 1 provides an algebraic tool to diagonalise the quaternion covariance matrix.
Observation 2 enables the diagonalisation of quaternion complementary covariance matrices, and degenerates into the Takagi factorisation of complex symmetric matrices if the quaternion vector $\mathbf{x}$ has two vanishing imaginary parts \cite{Zhang2006}.

\subsection{Diagonalisation of symmetric matrices}
The tools for the decomposition of quaternion pseudo-covariance matrices are still in their infancy, because:
\begin{itemize}
\item The tools for the analysis of the symmetric matrices in $\mathbb{C}$ cannot be readily generalised to $\mathbb{H}$. 
\item It is required to simultaneously diagonalise all the three complementary covariance matrices, a task with prohibitively many degrees of freedom. 
\end{itemize}
We shall start our analysis with the diagonalisation of $2 \times 2$ symmetric quaternion matrices.
\\\textbf{Proposition 1:} \emph{A $2 \times 2$ symmetric quaternion matrix
$\mathbf{A}$ admits the factorisation $\mathbf{A}=\mathbf{U}\boldsymbol{\Lambda}\mathbf{U}^T$ if either its diagonal elements or off-diagonal elements are real-valued},
where $\mathbf{U}$ is a quaternion unitary matrix and $\boldsymbol{\Lambda}$ is a real diagonal matrix with the singular values of $\mathbf{A}$ on the diagonal.
\begin{proof}
Based on the SVD of quaternion matrices \cite{Zhang1997}, $\mathbf{A}=\mathbf{U}\boldsymbol{\Lambda}\mathbf{V}^H$, the matrix
products $\mathbf{AA}^*$ and $\mathbf{A}^*\mathbf{A}$ can be expressed as
\begin{equation*}
\begin{split}
\mathbf{AA}^*=\mathbf{AA}^{H}=(\mathbf{U}\boldsymbol{\Lambda}\mathbf{V}^H) (\mathbf{U}\boldsymbol{\Lambda}\mathbf{V}^H)^{H}=\mathbf{U}\boldsymbol{\Lambda}^2 \mathbf{U}^{H}
\\
\mathbf{A}^*\mathbf{A}=\mathbf{A}^{H}\mathbf{A}=(\mathbf{U}\boldsymbol{\Lambda}\mathbf{V}^H)^H (\mathbf{U}\boldsymbol{\Lambda}\mathbf{V}^H)=\mathbf{V}\boldsymbol{\Lambda}^2 \mathbf{V}^{H}
\end{split}
\end{equation*}
The result in Appendix C shows that $(\mathbf{AA^*})^*=\mathbf{A}^*\mathbf{A}$ if either the diagonal or off-diagonal elements of $\mathbf{A}$ are real-valued. It then follows that $\mathbf{U}^*\boldsymbol{\Lambda}^2 (\mathbf{U}^*)^{H}=\mathbf{V}\boldsymbol{\Lambda}^2 \mathbf{V}^{H}$, where the commutativity of the product is valid because of the condition that the diagonal or off-diagonal elements of $\mathbf{A}$ are real-valued. Thus, we can obtain $\mathbf{V}=\mathbf{U}^*$, 
whereby the Takagi factorisation\footnote{The Takagi factorisation is a special case of the complex/quaternion SVD when $\mathbf{A}=\mathbf{A}^T$.} of quaternion symmetric matrices is identical to that of complex symmetric matrices, $\mathbf{A}=\mathbf{U}\boldsymbol{\Lambda}\mathbf{U}^T$.
\end{proof}
A consequence of Proposition 1 is that the Takagi factorisation of a quaternion symmetric matrix $\mathbf{A}$ is possible if $(\mathbf{AA^*})^*=\mathbf{A}^*\mathbf{A}$ holds. However, in general, $(\mathbf{AA^*})^*\neq\mathbf{A}^*\mathbf{A}$, indicating that such a Takagi factorisation does not exist for general quaternion symmetric matrices.
%
\subsection{Simultaneous diagonalisation of two covariance matrices} \label{sec Q-SUT}
We next proceed to introduce a simultaneous diagonalisation of the covariance and complementary covariance matrices. The finding that a non-singular transform is sufficient to diagonalise\footnote{Corollary 4.5.18 (c) in \cite{horn1990matrix} is the corresponding result in the complex domain.} both $\mathbf{C}_\mathbf{x}$ and $\mathbf{C}_{\mathbf{x}^{\alpha}}$ will form a basis for the quaternion uncorrelating transform~\cite{ShirinICASSP}. First, we shall clarify that the unitary transform which can simultaneously diagonalise $\mathbf{C}_\mathbf{x}$ and $\mathbf{C}_{\mathbf{x}^{\alpha}}$ for a general quaternion random vector $\mathbf{x}$ does not exist.
\\\textbf{Proposition 2:} \emph{If $\mathbf{C}_\mathbf{x}\mathbf{C}_{\mathbf{x}^{\alpha}}\neq(\mathbf{C}_\mathbf{x}\mathbf{C}_{\mathbf{x}^{\alpha}})^{\alpha H}$, then there is no quaternion unitary matrix $\mathbf{U}$ such that both $\mathbf{U}^H\mathbf{C}_\mathbf{x}\mathbf{U}$ and $\mathbf{U}^{H}\mathbf{C}_{\mathbf{x}^{\alpha}}\mathbf{U}^\alpha$ are diagonal.}
\begin{proof}
Define $\mathbf{\Lambda}_\mathbf{x}=\mathbf{U}^H\mathbf{C}_\mathbf{x}\mathbf{U}$, $\mathbf{\Lambda}_\alpha=\mathbf{U}^{H}\mathbf{C}_{\mathbf{x}^{\alpha}}\mathbf{U}^\alpha$, where $\mathbf{U}$ is a quaternion unitary matrix. It is then obvious that $\mathbf{\Lambda}_\mathbf{x}$ is a Hermitian matrix, so that if $\mathbf{\Lambda}_\mathbf{x}$ is diagonal, then $\mathbf{\Lambda}_\mathbf{x}$ is real diagonal. Furthermore, if $\mathbf{\Lambda}_\alpha$ is also diagonal, we can obtain
\begin{equation}\label{1}
\mathbf{C}_\mathbf{x}\mathbf{C}_{\mathbf{x}^{\alpha}}\!=\!\mathbf{U}\!\boldsymbol{\Lambda}_\mathbf{x}\!\mathbf{U}^H\!\mathbf{U}\!\boldsymbol{\Lambda}_\alpha\!\mathbf{U}^{\alpha H}\!\!=\!\!\mathbf{U}\!\boldsymbol{\Lambda}_\alpha\!\mathbf{U}^{\alpha H}\!\mathbf{U}^{\alpha}\!\boldsymbol{\Lambda}_\mathbf{x}\!\mathbf{U}^{\alpha H}\!=\!\mathbf{C}_{\mathbf{x}^{\alpha}}\mathbf{C}_\mathbf{x}^\alpha
\end{equation}
Because $\mathbf{C}_{\mathbf{x}}$ is Hermitian and $\mathbf{C}_{\mathbf{x}^{\alpha}}$ is $\alpha$-Hermitian, we can readily obtain
\begin{equation}\label{2}
\mathbf{C}_{\mathbf{x}^{\alpha}}\mathbf{C}_\mathbf{x}^\alpha=\mathbf{C}_{\mathbf{x}^{\alpha}}^{\alpha H}\mathbf{C}_\mathbf{x}^{\alpha H}=(\mathbf{C}_\mathbf{x}\mathbf{C}_{\mathbf{x}^{\alpha}})^{\alpha H}
\end{equation}
Combining \eqref{1} and \eqref{2} yields $\mathbf{C}_{\mathbf{x}}\mathbf{C}_{\mathbf{x}^\alpha}=(\mathbf{C}_\mathbf{x}\mathbf{C}_{\mathbf{x}^{\alpha}})^{\alpha H}$.
\end{proof}
Although a unitary transform for the simultaneous diagonalisation of $\mathbf{C}_\mathbf{x}$ and $\mathbf{C}_{\mathbf{x}^{\alpha}}$ is inapplicable, a non-unitary transform can be derived in a similar manner to the SUT in the complex domain. On the basis of Observation 1, $\mathbf{C}_\mathbf{x}$ can be factorised as
$
\mathbf{C}_\mathbf{x}=\mathbf{V}\mathbf{\Lambda}_\mathbf{x}\mathbf{V}^H
$
where $\mathbf{V}$ is a quaternion unitary matrix, and $\mathbf{\Lambda}_\mathbf{x}$ is a real diagonal matrix. Define a whitening transform $\mathbf{D}=\mathbf{V}\mathbf{\Lambda}_\mathbf{x}^{-\frac{1}{2}}\mathbf{V}^H$, and denote $\mathbf{s}=\mathbf{Dx}$, to obtain the covariance matrix of $\s$ as
\begin{equation*}
\mathbf{C}_\mathbf{s}=\mathbf{D}\mathbf{C}_\mathbf{x}\mathbf{D}^H=\mathbf{V}\mathbf{\Lambda}_\mathbf{x}^{-\frac{1}{2}}\mathbf{V}^H\mathbf{V}\mathbf{\Lambda}_\mathbf{x}\mathbf{V}^H\mathbf{V}\mathbf{\Lambda}_\mathbf{x}^{-\frac{1}{2}}\mathbf{V}^H=\mathbf{I}
\end{equation*}
Based on Observation 2, the $\alpha$-complementary covariance matrix of $\s$ can be factorised as $\mathbf{C}_{\mathbf{s}^\alpha}=\mathbf{W}\mathbf{\Lambda}_{\alpha}\mathbf{W}^{\alpha H}$,
where $\mathbf{W}$ is a quaternion unitary matrix, and $\mathbf{\Lambda}_{\alpha}$ is a real diagonal matrix. Now, the non-singular uncorrelating transform, $\mathbf{Q}=\mathbf{W}^H\mathbf{D}$, and its application in the form $\mathbf{y}=\mathbf{Qx}$ yield the simultaneous diagonalisation
\begin{align*}
\begin{split}
&\mathbf{C}_\mathbf{y}=\mathbf{W}^H\mathbf{C}_\mathbf{s}\mathbf{W}=
\mathbf{W}^H\mathbf{I}\mathbf{W}=\mathbf{I} \\
&\mathbf{C}_{\mathbf{y}^\alpha}\!=\!\!\mathbf{W}^H\mathbf{C}_{\mathbf{s}^\alpha}\mathbf{W}^\alpha
\!=\!\mathbf{W}^H\!\mathbf{W}\mathbf{\Lambda}_{\alpha}\mathbf{W}^{\alpha H}\!\mathbf{W}^\alpha\!=\!\mathbf{\Lambda}_{\alpha}
\end{split}
\end{align*}
We refer to the transform $\mathbf{Q}$ as the \emph{quaternion uncorrelating transform} (QUT), which is summarised in Proposition 3 and Algorithm~\ref{algorithm1}.
\\\textbf{Proposition 3 (QUT):}
\emph{For a random quaternion vector $\mathbf{x}$ with finite second-order statistics, there exists a quaternion non-singular matrix $\mathbf{Q}$ for which the transformed vector $\mathbf{y}=\mathbf{Qx}$ has the following covariance and $\alpha$-complementary covariance matrices:\begin{center}
$\mathbf{C}_{\mathbf{y}}=\mathbf{I}$  \hspace{5mm} $\mathbf{C}_{\mathbf{y}^\alpha}=\mathbf{\Lambda}_{\alpha}, \quad \text{where} \quad \alpha \in \{\imath,\jmath,\kappa\}$.\end{center}}
\begin{algorithm}
\begin{enumerate}
\item{Compute the eigendecomposition of the covariance matrix \newline $\mathbf{C}_\mathbf{x}=E\{\mathbf{x}\mathbf{x}^H\}=\mathbf{V} \boldsymbol{\Lambda}_\mathbf{x} \mathbf{V}^H$.}
\item{Compute the whitening matrix $\mathbf{D}=\mathbf{V}\mathbf{\Lambda}_\mathbf{x}^{-\frac{1}{2}}{\mathbf{V}}^H$.}
\item{Calculate the $\alpha$-complementary covariance of the whitened data, $\mathbf{s}=\mathbf{Dx}$, as $\mathbf{C}_{\mathbf{s}^\alpha}= E\{\mathbf{s}\mathbf{s}^{\alpha H}\}$}.
\item{Compute the factorisation of the $\alpha$-complementary covariance matrix, $\mathbf{C}_{\mathbf{s}^\alpha}=\mathbf{W}{\mathbf{\Lambda}_\alpha}\mathbf{W}^{\alpha H}$, where $\mathbf{W}$ is a quaternion unitary matrix and $\mathbf{\Lambda}_\alpha$ a real diagonal matrix.}
\item{The QUT matrix is then $\mathbf{Q}=\mathbf{W}^H\mathbf{D}$.}
\end{enumerate}
\caption{Quaternion Uncorrelating Transform (QUT).} \label{algorithm1}
\end{algorithm}
Recall that for a $\mathbb{C}^\alpha$-improper vector $\mathbf{x}$, only the $\alpha$-complementary covariance is non-vanishing while the other two complementary covariances vanish. Thus, the QUT of $\mathbb{C}^\alpha$-improper signals, for example, for the case $\alpha=\kappa$, satisfies
\begin{align*}
\mathbf{C}_\mathbf{y}=\mathbf{I} \hspace{7mm}
\mathbf{C}_{\mathbf{y}^\imath}=\mathbf{\Lambda}_\imath \hspace{7mm}
\mathbf{C}_{\mathbf{y}^\jmath}=\mathbf{\Lambda}_\jmath \hspace{7mm}
\mathbf{C}_{\mathbf{y}^\kappa}=\mathbf{\Lambda}_\kappa 
\end{align*}
\noindent where $\mathbf{\Lambda}_\kappa$ is a diagonal matrix and $\mathbf{\Lambda}_\imath=\mathbf{\Lambda}_\jmath=\mathbf{0}$. In other words, for a $\mathbb{C}^\alpha$-improper quaternion vector, the QUT can be regarded as the quaternion version of SUT \cite{de2002,Eriksson2006}. For a general improper quaternion vector, however, the remaining two complementary covariance matrices are still not diagonalised by the QUT. In fact, so far there is no closed-form solution to joint diagonalisation of four covariance matrices for general improper quaternion vectors. To circumvent this problem, we shall next introduce an approximate way to simultaneously diagonalise all four covariance matrices, a practically useful tool in most applications of quaternions, considering the effectiveness of the AUT in complex-valued signal processing \cite{mandic2015mean,xia2017complementary,xia2017full}.

\section{Simultaneous diagonalisation of all four covariance matrices}\label{sect:QAUT}
Following the approach in Section \ref{sect: diagonalisation}, we next propose the quaternion analogue of our proposed approximate uncorrelating transform (AUT) for complex-valued data \cite{CheongTook2012}. This will enable the simultaneous diagonlistion of the covariance matrix and three complementary covariance matrices.
\subsection{Univariate QAUT}
Consider first the univariate quaternion random vector $\mathbf{x}=\left(x_1, \ldots, x_L\right)^T$. Our aim is to find a unitary transform
$\mathbf{y}=\mathbf{\Phi}\mathbf{x}$ which simultaneously diagonalises the covariance matrix, $\mathbf{C}_{\mathbf{y}}=E\left\{ \mathbf{y}\mathbf{y}^{H}\right\}=\mathbf{\Phi}\mathbf{C}_{\mathbf{x}}\mathbf{\Phi}^H$, and the three complementary covariance matrices, $\mathbf{C}_{\mathbf{y}^\alpha}=E\left\{ \mathbf{y}\mathbf{y}^{\alpha H}\right\}=\mathbf{\Phi}\mathbf{C}_{\mathbf{x}^\alpha}\mathbf{\Phi}^{\alpha H}$, $\alpha \in \{\imath,\jmath,\kappa\}$.
Notice that if $\mathbf{C}_{\mathbf{x}}\mathbf{C}_{\mathbf{x}}^H=\mathbf{C}_{\mathbf{x}^\alpha}\mathbf{C}_{\mathbf{x}^\alpha}^H$, we have $\mathbf{C}_{\mathbf{y}}\mathbf{C}_{\mathbf{y}}^H=\mathbf{C}_{\mathbf{y}^\alpha}\mathbf{C}_{\mathbf{y}^\alpha}^H$.

For $N$ samples of $\mathbf{x}$, $N>L$, the covariance matrix, and the $\alpha$-complementary covariance matrix can be estimated as 
\begin{equation}\label{covariance}
\mathbf{C}_{\mathbf{x}}=\frac{1}{N}\!\underset{n_1=1}{\overset{N}{\sum}}\mathbf{x}\left(n_1\right)\mathbf{x}^{H}\left(n_1\right)
\end{equation}
\begin{equation}\label{complementary covariance}
\mathbf{C}_{\mathbf{x}^\alpha}=\frac{1}{N}\!\underset{n_1=1}{\overset{N}{\sum}}\mathbf{x}\left(n_1\right)\mathbf{x}^{\alpha H}\left(n_1\right)
\end{equation}
Their quadratic forms $\mathbf{C}_{\mathbf{x}}\mathbf{C}_{\mathbf{x}}^H$ and $\mathbf{C}_{\mathbf{x}^\alpha}\mathbf{C}_{\mathbf{x}^\alpha}^H$ are calculated as
\small
\begin{align*}
&\mathbf{C}_{\mathbf{x}}\mathbf{C}_{\mathbf{x}}^{H} \\\!\!\!&\!\!=\!\!\left[\!\frac{1}{N}\!\underset{n_1=1}{\overset{N}{\sum}}\mathbf{x}\left(n_1\right)\mathbf{x}^{H}\left(n_1\right)\!\right]\!\!\left[\!\frac{1}{N}\!\underset{n_2=1}{\overset{N}{\sum}}\mathbf{x}\left(n_2\right)\mathbf{x}^{H}\left(n_2\right)\!\right]\\
\!\!&\!\!=\!\!\frac{1}{N^{2}}\underset{n_1=1}{\overset{N}{\sum}}\underset{n_2=1}{\overset{N}{\sum}}\mathbf{x}\left(n_1\right)\!\!\left[\underset{l=1}{\overset{L}{\sum}}x_{l}^{*}\left(n_1\right)x_{l}\left(n_2\right)\right]\!\!\mathbf{x}^{H}\left(n_2\right)
\end{align*}
\begin{align*}
&\mathbf{\mathbf{C}}_{\mathbf{x}^\alpha}\mathbf{C}_{\mathbf{x}^\alpha}^{H}\\
&\!\!=\!\!\left[\!\frac{1}{N}\!\!\underset{n_1=1}{\overset{N}{\sum}}\mathbf{x}\left(n_1\right)\mathbf{x}^{\alpha H}\left(n_1\right)\!\right]\!\!\left[\!\frac{1}{N}\!\!\underset{n_2=1}{\overset{N}{\sum}}\mathbf{x}^{\alpha}\left(n_2\right)\mathbf{x}^{H}\left(n_2\right)\!\right]\\
&\!\!=\!\!\frac{1}{N^{2}}\underset{n_1=1}{\overset{N}{\sum}}\underset{n_2=1}{\overset{N}{\sum}}\mathbf{x}\left(n_1\right)\!\!\left[\underset{l=1}{\overset{L}{\sum}}x_{l}^{*}\left(n_1\right)x_{l}\left(n_2\right)\right]^{\alpha}\!\!\!\!\mathbf{x}^{H}\left(n_2\right)
\end{align*}
\normalsize
Upon applying the approximation
\begin{equation}\label{single condition}
\underset{l=1}{\overset{L}{\sum}}x_{l}^{*}\left(n_1\right)x_{l}\left(n_2\right)
\!\!\approx\!\!\left[\underset{l=1}{\overset{L}{\sum}}x_{l}^{*}\left(n_1\right)x_{l}\left(n_2\right)\right]^{\alpha}
\end{equation}
we obtain $\mathbf{C}_{\mathbf{x}}\mathbf{C}_{\mathbf{x}}^H\approx\mathbf{C}_{\mathbf{x}^\alpha}\mathbf{C}_{\mathbf{x}^\alpha}^H$, and hence $\mathbf{C}_{\mathbf{y}}\mathbf{C}_{\mathbf{y}}^H\approx\mathbf{C}_{\mathbf{y}^\alpha}\mathbf{C}_{\mathbf{y}^\alpha}^H$, which implies that if $\mathbf{C}_{\mathbf{y}}$ is diagonal, so too is $\mathbf{C}_{\mathbf{y}^\alpha}$, and \emph{vice versa}. Furthermore, if 
\begin{equation}\label{approximation}
\begin{split}
\underset{l=1}{\overset{L}{\sum}}x_{l}^{*}(n_1)x_{l}(n_2)\in \mathbb{R}
\end{split}
\end{equation} 
holds approximately, then \eqref{single condition} holds for all $\alpha \in \{ \imath, \jmath, \kappa \}$, and thus the simultaneous diagonalisation applies to all of the three complementary covariance matrices, producing the following \emph{quaternion approximate uncorrelating transform} (QAUT). 
\\\textbf{Proposition 4 (QAUT):}
\emph{If any of $\mathbf{C}_{\mathbf{x}}$, $\mathbf{C}_{\mathbf{x}^\imath}$,
$\mathbf{C}_{\mathbf{x}^\jmath}$ and $\mathbf{C}_{\mathbf{x}^\kappa}$ is diagonalised by the unitary matrix $\mathbf{\Phi}$, the structural similarity between the four covariance matrices, displayed in \eqref{single condition},} makes possible the approximate diagonalisation of the other three covariance matrices by $\mathbf{\Phi}$. From \eqref{eq: covariance relationship}, the pseudo-covariance matrix is also approximately diagonalised in this way.

The transformation matrix $\mathbf{\Phi}$ can take the form of $\mathbf{Q}$ in Observation 1 or $\mathbf{Q}_\alpha$ in Observation 2. 
For example, the implementation of the transform $\mathbf{Q}$, obtained from the eigendecomposition of $\mathbf{C}_\mathbf{x}$, yields
\begin{equation}\label{QAUT}
\begin{split}
\mathbf{C}_\mathbf{x}&=\mathbf{Q}^{H}\mathbf{\boldsymbol{\Lambda}}_{\mathbf{x}}\mathbf{Q}\\
\mathbf{C}_\mathbf{x^\imath}&\approx\mathbf{Q}^{H}\mathbf{\boldsymbol{\Lambda}}_{\mathbf{\imath}}\mathbf{Q}^{\imath}\\
\mathbf{C}_\mathbf{x^\jmath}&\approx\mathbf{Q}^{H}\mathbf{\boldsymbol{\Lambda}}_{\mathbf{\jmath}}\mathbf{Q}^{\jmath}\\
\mathbf{C}_\mathbf{x^\kappa}&\approx\mathbf{Q}^{H}\mathbf{\boldsymbol{\Lambda}}_{\mathbf{\kappa}}\mathbf{Q}^{\kappa}
\end{split}
\end{equation}
\subsection{Principle of QAUT}
The QAUT is based on the approximation in \eqref{approximation}. To examine the validity of this approximation for highly improper quaternion signals, assume that the four components of the quaternion variable $x_l$ are perfectly correlated and consider the Cayley-Dickson construction \cite{Ward1997} \[x_l=x_{a}+x_{b} \imath + x_{c} \jmath+x_{d} \kappa =c_{1}+c_{2}\jmath\]
where $c_{1}=x_{a}+x_{b} \imath $ and $c_{2}=x_{c}+x_{d} \imath $ are complex variables defined in the plane spanned by $\{1,\imath\}$. This gives
\begin{align}
\begin{array}{rl}
\!\!\!\!\!x^{*}_l\left(n_1\right)x_l\left(n_2\right)=\!\!&\!\!c_{1}^{*}\left(n_1\right)c_{1}\left(n_2\right)+c_{2}\left(n_1\right)c_{2}^{*}\left(n_2\right)+\\&\!\!\left[c_{1}^{*}\left(n_1\right)c_{2}\left(n_2\right)-c_{2}\left(n_1\right)c_{1}^{*}\left(n_2\right)\right]\jmath
\end{array}\label{x_k*x_l-original}
\end{align}
Upon applying the statistical expectation operator, we have
\begin{align}
\begin{array}{rl}
\!\!\!\!\!\!E\left\{x^{*}_l\left(n_1\right)x_l\left(n_2\right)\right\}=&E\left\{ c_{1}^{*}\left(n_1\right)c_{1}\left(n_2\right)\right\}+E\left\{ c_{2}\left(n_1\right)c_{2}^{*}\left(n_2\right)\right\}\\&+E\left\{ c_{1}^{*}\left(n_1\right)c_{2}\left(n_2\right)-c_{2}\left(n_1\right)c_{1}^{*}\left(n_2\right)\right\} \jmath 
\end{array}\label{x_k*x_l}
\end{align}
where
\begin{align*}
\!\!E\{c_{1}^{*}\left(n_1\right)c_{1}\left(n_2\right)\}=&E\{x_a\left(n_1\right)x_a\left(n_2\right)+x_b\left(n_1\right)x_b\left(n_2\right)\}+\\
&E\{x_a\left(n_1\right)x_b\left(n_2\right)-x_a\left(n_2\right)x_b\left(n_1\right)\}\imath\\
\approx&E\{x_a\left(n_1\right)x_a\left(n_2\right)+x_b\left(n_1\right)x_b\left(n_2\right)\}+\\
&E\{x_b\left(n_1\right)x_b\left(n_2\right)-x_b\left(n_2\right)x_b\left(n_1\right)\}\imath\\
=&E\{x_a\left(n_1\right)x_a\left(n_2\right)+x_b\left(n_1\right)x_b\left(n_2\right)\}\\\in & \: \mathbb{R}\numberthis
\end{align*}
and likewise for $E\{c_{2}\left(n_1\right)c_{2}^*\left(n_2\right)\}\in\mathbb{R}$. 
Moreover, if $c_{1}$ and $c_{2}$ are highly correlated, they can be expressed as $c_{2}\left(t\right)\approx kc_{1}\left(t\right)+g$, where $k$ and $g$ are complex-valued constants, to give
\begin{align*}
& E\left\{ c_{1}^{*}\left(n_1\right)c_{2}\left(n_2\right)-c_{2}\left(n_1\right)c_{1}^{*}\left(n_2\right)\right\}  \\
&\approx E\left\{ c_{1}^{*}\left(n_1\right)\left[kc_{1}\left(n_2\right)+g\right]-\left[kc_{1}\left(n_1\right)+g\right]c_{1}^{*}\left(n_2\right)\right\} \\
&= 2kE\left\{\mathfrak{I}\left[c_{1}^{*}\left(n_1\right)c_{1}\left(n_2\right)\right]\right\}+E\{g\}E\left\{c_{1}^{*}\left(n_1\right)\right\}\\&\quad-E\{g\}E\left\{c_{1}^{*}\left(n_2\right)\right\}\\&\approx0 \numberthis \label{11}
\end{align*}
Therefore, we obtain approximately
\begin{equation}
E\left\{x^{*}_l\left(n_1\right)x_l\left(n_2\right)\right\}\in\mathbb{R} \label{real}
\end{equation} 
and the approximation \eqref{approximation} holds given the assumption 
\begin{equation*}
E\left\{x^{*}_l\left(n_1\right)x_l\left(n_2\right)\right\}\approx \frac{1}{L} \underset{l=1}{\overset{L}{\sum}}x_{l}^{*}(n_1)x_{l}(n_2)
\end{equation*} 
%

From the above analysis, the QAUT holds exactly when
the four components
of the quaternion variable are perfectly correlated (a maximally improper $x$). For general quaternion signals, however, these components are not perfectly correlated, causing the approximation error impacted by the correlation coefficients, given by
\begin{align*}
\rho_{ab}=\frac{\textrm{cov}[x_{a},x_{b}]}{\sigma_{x_{a}}\sigma_{x_{b}}},\quad\rho_{ac}=\frac{\textrm{cov}[x_{a},x_{c}]}{\sigma_{x_{a}}\sigma_{x_{c}}},\quad\rho_{ad}=\frac{\textrm{cov}[x_{a},x_{d}]}{\sigma_{x_{a}}\sigma_{x_{d}}},
\\\!\!\!\!\!\quad\rho_{bc}=\frac{\textrm{cov}[x_{b},x_{c}]}{\sigma_{x_{b}}\sigma_{x_{c}}},\quad\rho_{bd}=\frac{\textrm{cov}[x_{b},x_{d}]}{\sigma_{x_{b}}\sigma_{x_{d}}},\quad\rho_{cd}=\frac{\textrm{cov}[x_{c},x_{d}]}{\sigma_{x_{c}}\sigma_{x_{d}}}
\end{align*}
with the range $[-1,1]$, where `$\textrm{cov}$' denotes the covariance and $\sigma$ the standard deviation.
From the diagonalisation condition in \eqref{single condition}, the diagonalisation error of $\C_{\x^\alpha}$ arises from and is positively correlated with
\begin{equation}\label{error}
\begin{split}
&x^{*}_l\left(n_1\right)x_l\left(n_2\right)-\left[x^{*}_l\left(n_1\right)x_l\left(n_2\right)\right]^\alpha
\\&=2\mathfrak{I}_\beta\left[x_l^{*}\left(n_1\right)x_l\left(n_2\right)\right]\beta+2\mathfrak{I}_\gamma\left[x_l^{*}\left(n_1\right)x_l\left(n_2\right)\right]\gamma
\\&\text{for distinct}~~\alpha,\beta,\gamma \in \{\imath,\jmath,\kappa\}\end{split}
\end{equation}
Note that the imaginary part of $x_l^{*}\left(n_1\right)x_l\left(n_2\right)$ is given by
\begin{align}\label{Im}
\begin{array}{rl}
&\mathfrak{I}\left[x_l^{*}\left(n_1\right)x_l\left(n_2\right)\right]\\&=\!\!\left[x_{a}\left(n_1\right)x_{b}\left(n_2\right)-x_{a}\left(n_2\right)x_{b}\left(n_1\right)+x_{c}\left(n_2\right)x_{d}\left(n_1\right)-x_{c}\left(n_1\right)x_{d}\left(n_2\right)\right]\imath\\&+\left[x_{a}\left(n_1\right)x_{c}\left(n_2\right)-x_{a}\left(n_2\right)x_{c}\left(n_1\right)+x_{b}\left(n_1\right)x_{d}\left(n_2\right)-x_{b}\left(n_2\right)x_{d}\left(n_1\right)\right]\jmath\\&+\left[x_{a}\left(n_1\right)x_{d}\left(n_2\right)-x_{a}\left(n_2\right)x_{d}\left(n_1\right)+x_{b}\left(n_2\right)x_{c}\left(n_1\right)-x_{b}\left(n_1\right)x_{c}\left(n_2\right)\right]\kappa
\end{array}
\end{align}
Similarly to the AUT for complex-valued data \cite{CheongTook2012}, equation \eqref{Im} can be rewritten as 
\begin{align}\label{Im_1}
\!\!\!\!\begin{array}{rl}
\mathfrak{I}\!\left[x_l^{*}\left(n_1\right)x_l\left(n_2\right)\right]\!=\!&\!\!\!\!\left(\!\xi_{ab}\sqrt{1-\rho_{ab}^{2}}+\xi_{cd}\sqrt{1-\rho_{cd}^{2}}\!\right)\!\!\imath\!+\!\left(\!\xi_{ac}\sqrt{1-\rho_{ac}^{2}}+\xi_{bd}\sqrt{1-\rho_{bd}^{2}}\!\right)\!\!\jmath\!\\&\!\!\!\!+\left(\xi_{ad}\sqrt{1-\rho_{ad}^{2}}+\xi_{bc}\sqrt{1-\rho_{bc}^{2}}\right)\kappa
\end{array}
\end{align}
where $\xi_{ab},\xi_{cd},\xi_{ac},\xi_{bd},\xi_{ad},\xi_{bc}$ are real-valued coefficients dependent on the statistics of $x_l$. 

It is now obvious that when the components of quaternion data become more correlated, \eqref{Im_1} decreases and the diagonalisation error hence decreases.
Owing to the cumulative nature of the error in \eqref{approximation}, the approximation error also
increases with the data segment length $L$.
\subsection{Multivariate QAUT}
We next extend the univariate QAUT to the multivariate case when an $M$-variate quaternion signal is represented by a
data matrix
$\mathbf{X}=\left[\begin{array}{ccc}
\mathbf{x}_{1},\ldots,\mathbf{x}_{M}\end{array}\right]^{T}$, where the $m$-th column 
random vector $\mathbf{x}_m=\left[\begin{array}{ccc}x_{m}(1),\ldots,x_{m}(N)\end{array}\right]^{T}$ ($N>M\geq m$) represents $N$ samples of the $m$-th variate. The covariance and $\alpha$-covariance matrices of the $M$ variates are then given by $\mathbf{C}_\mathbf{X}=\frac{1}{N} E\{\mathbf{X}\mathbf{X}^{H}\}$ and
$\mathbf{C}_\mathbf{X^\alpha}=\frac{1}{N} E\{\mathbf{X}\mathbf{X}^{\alpha H}\}$. Our aim is to find a unitary transform
\begin{equation}
\mathbf{Y}=\mathbf{\Phi}\mathbf{X} \label{M-QAUT}
\end{equation}
which simultaneously diagonalises the covariance matrix, $\mathbf{C}_{\mathbf{Y}}=E\{\mathbf{Y}\mathbf{Y}^{H}\}=\mathbf{\Phi}\mathbf{C}_{\mathbf{X}}\mathbf{\Phi}^{H}$, and the three complementary covariance matrices, $\mathbf{C}_{\mathbf{Y}^\alpha}=E\{\mathbf{Y}\mathbf{Y}^{\alpha H}\}=\mathbf{\Phi}\mathbf{C}_{\mathbf{X}^\alpha}\mathbf{\Phi}^{\alpha H}$, $\alpha \in \{\imath,\jmath,\kappa\}$. Notice that if $\mathbf{C}_{\mathbf{X}}\mathbf{C}_{\mathbf{X}}^H=\mathbf{C}_{\mathbf{X}^\alpha}\mathbf{C}_{\mathbf{X}^\alpha}^H$, then $\mathbf{C}_{\mathbf{Y}}\mathbf{C}_{\mathbf{Y}}^H=\mathbf{C}_{\mathbf{Y}^\alpha}\mathbf{C}_{\mathbf{Y}^\alpha}^H$.  

The covariance and $\alpha$-complementary covariance matrices can be estimated from $\X$ as 
\begin{equation}\label{covariance-1}
\mathbf{C}_\mathbf{X}=\frac{1}{N}\mathbf{X}\mathbf{X}^{H}
\end{equation}
\begin{equation}\label{complementary covariance-1}
\mathbf{C}_\mathbf{X^\alpha}=\frac{1}{N}\mathbf{X}\mathbf{X}^{\alpha H}
\end{equation}
Their quadratic forms $\mathbf{C}_{\mathbf{X}}\mathbf{C}_{\mathbf{X}}^H$ and $\mathbf{C}_{\mathbf{X}^\alpha}\mathbf{C}_{\mathbf{X}^\alpha}^H$ are calculated as
\begin{align*}
\begin{array}{rl}
\mathbf{C}_\mathbf{X}\mathbf{C}_\mathbf{X}^{H}=&\!\!\frac{1}{N^2}\mathbf{X}\mathbf{X}^{H}\mathbf{X}\mathbf{X}^{H}=\frac{1}{N^2}\mathbf{X}(\underset{m=1}{\overset{M}{\sum}}\mathbf{x}_{m}^{*}\mathbf{x}_{m}^T)\mathbf{X}^{H}\\
\mathbf{C}_\mathbf{X^\alpha}\mathbf{C}_\mathbf{X^\alpha}^{H}=&\!\!\frac{1}{N^2}\mathbf{X}\mathbf{X}^{\alpha H}\mathbf{X}^{\alpha}\mathbf{X}^{H}=\frac{1}{N^2}\mathbf{X}(\underset{m=1}{\overset{M}{\sum}}\mathbf{x}_{m}^{*}\mathbf{x}_{m}^T)^\alpha\mathbf{X}^{H}
\end{array}
\end{align*}
Following the analysis in \eqref{x_k*x_l-original}-\eqref{11}, we can prove that the $[n_1,n_2]$ element of the matrix $E\{\mathbf{x}_m^{*}\mathbf{x}_m^T\}$, $E\{x_{m}^{*}(n_1)x_{m} (n_2)\}$, is real-valued, for all $n_1,n_2 \in \{1,\ldots,N\}$, whereby $E\{\mathbf{x}_m^{*}\mathbf{x}_m^T\}\in \mathbb{R}^{N\times N}$ and  
\begin{equation*}
\underset{m=1}{\overset{M}{\sum}}\!\mathbf{x}_{m}^{*}\mathbf{x}^T_{m}\in \mathbb{R}^{N\times N}
\end{equation*}
Under this condition, we obtain $\mathbf{C}_{\mathbf{x}}\mathbf{C}_{\mathbf{x}}^H\approx\mathbf{C}_{\mathbf{x}^\alpha}\mathbf{C}_{\mathbf{x}^\alpha}^H$, so that $\mathbf{C}_{\mathbf{Y}}\mathbf{C}_{\mathbf{Y}}^H\approx\mathbf{C}_{\mathbf{Y}^\alpha}\mathbf{C}_{\mathbf{Y}^\alpha}^H$, which implies if $\mathbf{C}_{\mathbf{Y}}$ is diagonal, so too is $\mathbf{C}_{\mathbf{Y}^\alpha}$, and \emph{vice versa}. This means that the four covariance matrices of multivariate quaternion data can be simultaneously diagonalised via the QAUT, such as in \eqref{QAUT}.

Similarly to the univariate case, the diagonalisation error of the multivariate QAUT decreases with an increase in the correlation between the components of quaternion data, and increases with the number of variates, $M$.
\begin{rem} 
The QAUT becomes exact for quaternion signals with special second-order statistical properties, such as $\mathbb{H}$-proper signals, the signals with fully correlated components, and the signals  
for which the covariance and complementary covariance matrices have the same singular vectors. This can be verified by exploring the underlying matrix structures.
\end{rem}
\section{Simulation studies}\label{sect: simulation}
The effectiveness of the proposed QUT and QAUT techniques is illustrated via simulations on both synthetic and real-world quaternion signals. In the experiments, the covariance and complementary covariance matrices were estimated based on \eqref{covariance}, \eqref{complementary covariance}, \eqref{covariance-1} and \eqref{complementary covariance-1}.  
\subsection{Performance of QUT for $\mathbb{C}^\eta$-improper signals}
In the first experiment, the QUT was applied to the decorrelation of multivariate $\mathbb{C}^\kappa$-improper quaternion signals. First, three uncorrelated $\mathbb{C}^\kappa$-improper quaternion signals were generated through the approach introduced in Appendix B and were subsequently mixed
through a $3\times3$ matrix for which the elements were drawn from a standard normal
distribution. The mixed $\mathbb{C}^\kappa$-improper signals, $\mathbf{x}_1$, $\mathbf{x}_2$ and $\mathbf{x}_3$, were correlated in terms of the $\kappa$-complementary covariance, with the correlation coefficients 0.69 (between $\mathbf{x}_1$ and $\mathbf{x}_2$), 0.34 (between $\mathbf{x}_1$ and $\mathbf{x}_3$), and 0.91 (between $\mathbf{x}_2$ and $\mathbf{x}_3$). Fig. 1 shows 3-D scatter plots of the mixed signals, $\mathbf{x}_1$, $\mathbf{x}_2$ and $\mathbf{x}_3$, illustrating their high correlation, as indicated by elliptical shapes of the scatter plots at an angle to the coordinate axes. The QUT was used to decorrelate $\mathbf{x}_1$, $\mathbf{x}_2$ and $\mathbf{x}_3$ into $\mathbf{y}_1$, $\mathbf{y}_2$ and $\mathbf{y}_3$, as illustrated in Fig. 2, which shows much more circular nature of $\mathbf{y}_1$, $\mathbf{y}_2$ and $\mathbf{y}_3$, as desired.
\subsection{Performance of QAUT for synthetic improper signals}
To assess the performance of QAUT, the squared diagonal error of the complementary covariance matrix $\mathbf{C}_{\mathbf{x}^\alpha}$ was measured by an average power ratio of the off-diagonal,
$\lambda_{\alpha,ij}$, versus diagonal, $\lambda_{\alpha,ii}$, elements of the
approximately diagonal matrix $\mathbf{\boldsymbol{\Lambda}}_{\alpha}$, given by
\[
\varepsilon_{\alpha}^{2}=\frac{\underset{i,j=1}{\overset{L}{\sum}}\left|\lambda_{\alpha,ij}\right|^{2}}{(N-1)\underset{i=1}{\overset{L}{\sum}}\lambda_{\alpha,ii}^{2}}\times100\%,\quad i \neq j
\]

Synthetic quaternion signals with a varying degree
of correlation between the real, $\imath$, $\jmath$ and $\kappa$ components were considered. For
simplicity, we assumed the six correlation coefficients of the four components to be equally distributed in $(0,1]$ and denoted by $\rho$, where $\rho=1$ indicates full correlation and $\rho=0$ no correlation.

In the univariate simulations, a synthetic quaternion data vector with a varying $\rho$ was generated from a quaternion white Gaussian
signal. Then, the performance of the univariate QAUT was assessed comprehensively against the degree of correlation present in the data components, $\rho\in\left[0.5,1\right]$, and the length of data segment, $L\in\left[5,20\right]$. Conforming with the analysis, Fig. \ref{fig:synthetic_univariate} shows that the diagonalisation error increased as $\rho$ or $L$ increased. The performance was excellent, for example, for a moderate correlation degree, $\rho=0.5$, and a long data segment, $L=20$, where the squared diagonalisation error was less than $1\%$, while for a high correlation degree, $\rho$ is close to unit, where the error was negligible. 

In the multivariate simulations, multiple channels of quaternion signals with a varying $\rho$ were generated and then mixed using a matrix the elements of which were drawn from a standard normal distribution, to generate multivariate quaternion data. The performance of the multivariate QAUT was assessed against the degree of correlation present in the data components, $\rho\in\left[0.5,1\right]$, and the number of variates, $M\in\left[5,20\right]$. Fig. \ref{fig:synthetic_multivariate} shows that the diagonalisation error increased as $\rho$ or $M$ increased. Similarly to the univariate case, for a moderate correlation degree, $\rho=0.5$, and a large number of data variates, $M=20$, the squared diagonalisation error was less than $1\%$, while for a high correlation degree, $\rho$ is close to unit, the error was negligible. 
\subsection{Performance of QAUT for real-world signals}
The QAUT was also applied to 4-D real-world electroencephalogram (EEG) signals by combining four adjacent channels of real-valued EEG signals into a quaternion-valued signal. The four channels measuring adjacent brain regions exhibited strong correlation with each other while the channels measuring distant brain regions were relatively weakly correlated \cite{Park2014,enshaeifar2016quaternion}. We tested the QAUT on such quaternion EEG signals with low, medium and high correlations between the four dimensions, with the data segment length $L$ varying between 5 and 20. Fig. \ref{fig:EEG} illustrates that the diagonalisation error increased with either an increase in $L$ or a decrease in the degree of correlation present in the data. The squared diagonal error was less than $0.2\%$.
\subsection{QAUT for rank reduction}
In many practical applications the acquired data are highly correlated or even collinear, with redundant elements that do not contribute to the accuracy of processing but require excessive computational resources and may affect stability of algorithms through the associated singular correlation matrices. It is therefore desirable to use the minimum number of variables to describe the information within a data set; this is typically achieved through a low-rank matrix approximation. 
For data requiring rank reduction, the assumption that the components of variables are highly correlated is sensible, that is, the approximation conditions for the QAUT hold. 
The multivariate QAUT in \eqref{M-QAUT} can be therefore rearranged to decompose the data matrix $\mathbf{X}$ into a sum of uncorrelated rank one vector outer products 
\[\mathbf{X}=\boldsymbol{\Phi}^H\mathbf{Y}=\underset{m=1}{\overset{M}{\sum}}\boldsymbol{\phi}_m^{H}\mathbf{y}_m\]
where $\boldsymbol{\phi}_m$ and $\mathbf{y}_m$ are rows in the matrices $\boldsymbol{\Phi}$ and $\mathbf{Y}$, respectively, and the transformed variables, $\mathbf{y}_m$, are monotonically non-increasing in the variance. The use of variables for which the variance is above a defined threshold and the disposal of the remaining variables, which account for noise, provides the reduced rank approximation of $\mathbf{X}$ in the form
\[\hat{\mathbf{X}}=\underset{m=1}{\overset{P}{\sum}}\boldsymbol{\phi}_m^{H}\mathbf{y}_m,  \quad \quad \text{where} \quad P<M\]
which is illustrated in Fig. \ref{fig:QAUTrankred}. This is analogous to the real domain where the SVD provides a change of coordinates to yield uncorrelated variables for which the variance is monotonically non-increasing.

To illustrate the potential of the QAUT in rank reduction of quaternion data, a 2-D quaternion process, $\mathbf{x}_1$, $\mathbf{x}_2$, was generated where $\mathbf{x}_1$ and $\mathbf{x}_2$ were highly correlated as shown in Fig. \ref{fig:RankRed1}, which plots the component relationship in the real, $\imath$, $\jmath$ and $\kappa$ planes. The QAUT decorrelated $\mathbf{x}_1$, $\mathbf{x}_2$ into $\mathbf{y}_1$ and $\mathbf{y}_2$. Fig. \ref{fig:RankRed2} shows significant variation only in the $\mathbf{y}_2$ axis while the $\mathbf{y}_1$ variation corresponds to noise, indicating that the true dimensionality of the process is one rather than two.

Rank reduction is commonly used in chemometrics where many sets of measurements are required to examine a process and the important information may not be directly observable. For example, Near Infrared Spectra (NIR) are used to examine the chemical content of a data sample and are often highly correlated or collinear. We applied the QAUT to the NIR obtained from 80 samples of corn using three different spectrometers \cite{Corndata}. 
A pure quaternion variable was formed for each corn sample by combining its spectra from the three spectrometers. The QAUT diagonalised the data covariance matrices successfully to show that one quaternion variable was sufficient to express most of the information, as the largest diagonal element was 39 compared to 0.01 for the next largest, as shown in Fig. \ref{fig:RankRedCorn}, which displays the logarithm of the ratio of variance explained by each variable, compared to the total variance in the data. The insert plot in Fig. \ref{fig:RankRedCorn} shows the cumulative ratio of variance explained by including further variables for the original and transformed data. Observe that the transformed data with much fewer variables explain the same information as the original, thus verifying the QAUT.
\subsection{QAUT for imbalance detection in three-phase power systems}
As shown in Section \ref{sect:QAUT}, when the components of data are highly correlated, the QAUT diagonalises all complementary covariance matrices at the expense of a small error, however, this error increases for lower degrees of correlation. We now show a way to exploit this data dependence of the QAUT and the structure of the complementary covariance matrices to detect imbalance in a three-phase power system. 

In general, the instantaneous voltages of a three-phase power system are given by \cite{Talebi2015a}
\begin{align}\label{eq:3phase}
\begin{array}{rl}
 v_{1}(n) &= A_1\sin(2\pi f\Delta T n + \varphi_1) \\ 
 v_{2}(n) &= A_2\sin(2\pi f\Delta T n + \varphi_2) \\ 
 v_{3}(n) &= A_3\sin(2\pi f\Delta T n + \varphi_3) 
 \end{array}
\end{align}
where $A_1, A_2, A_3$ are the instantaneous amplitudes, $\varphi_1, \varphi_2, \varphi_3$ the instantaneous phases, $f$ the system frequency, and $\Delta T$ the sampling interval. In a balanced three-phase system, $A_1=A_2=A_3$, $\varphi_2-\varphi_1=\frac{2\pi}{3}$, $\varphi_3-\varphi_2=\frac{2\pi}{3}$, and the correlation between each component is equivalent \cite{Xia2012}. 
The power grid is usually designed to operate optimally under balanced conditions; however, faults in the power system can cause imbalanced operating conditions that propagate through the network, threatening its stability. Therefore, it is important in fault detection and mitigation applications to identify the incidences when the power grid is operating in an unbalanced fashion \cite{Kanna2015}. 

Pure quaternions have been used to deal with data recordings from all three phases simultaneously as $x(n)=v_{1}(n) \imath + v_{2}(n) \jmath +v_{3}(n) \kappa $ \cite{Talebi2015a}. Here we organise the recorded data into a pure quaternion vector, $\x = \v_1 \imath + \v_2 \jmath +  \v_3 \kappa$, where $\v_m=[v_m(1),v_m(2),\ldots,v_m(n)]^T$, $m=1,2,3$. We can then apply the univariate QAUT to $\x$ and employ the diagonalisation errors of the complementary covariance matrices as an indicator of power imbalance. This can be explained by the overall diagonalisation error given by 
\begin{align}\label{Error_power}
\begin{array}{rl}
\!\!\!\!\!\!\!\mathfrak{I}\left[x_l^{*}\left(n_1\right)x_l\left(n_2\right)\right]\!\!=\!\!\!&\!\!\!\left[v_{2}\left(n_2\right)v_{3}\left(n_1\right)\!\!-\!\!v_{2}\left(n_1\right)v_{3}\left(n_2\right)\right]\imath\\\!\!&\!\!\!\!+\!\left[v_{1}\left(n_1\right)v_{3}\left(n_2\right)\!\!-\!\!v_{1}\left(n_2\right)v_{3}\left(n_1\right)\right]\!\jmath\\\!\!&\!\!\!\!+\!\left[v_{1}\left(n_2\right)v_{2}\left(n_1\right)\!\!-\!\!v_{1}\left(n_1\right)v_{2}\left(n_2\right)\right]\!\kappa
\end{array}
\end{align} 
which is a reduced form of \eqref{Im} and is positively correlated with the diagonalisation error. When the system is balanced, the covariances $E\{\v_p\v_q^T\}$ are equivalent for all phases $p, q \in \{1, 2, 3\}$, so the three imaginary parts of $x_l^{*}\left(n_1\right)x_l\left(n_2\right)$ in \eqref{Error_power} are equal in absolute value, leading to equal diagonalisation errors of the three complementary covariance matrices. On the other hand, when the system becomes unbalanced, the correlation degrees between the three phase voltages are different, and hence the three imaginary parts of $x_l^{*}\left(n_1\right)x_l\left(n_2\right)$ differ, resulting in different diagonalisation errors of the three complementary covariance matrices. This allows the QAUT to be used in an imbalance detection application. To this end, we applied the QAUT to a sliding window of 0.04 seconds of a three-phase power system for which a fault occurred at 0.25s lasting for 0.25s before the system returned to a balanced state. During this fault, $\v_1$ decreased significantly in amplitude whereas $\v_2$ and $\v_3$ increased in amplitude ($\v_2 > \v_3$). Fig. \ref{fig:powimbal} shows that the diagonalisation error for $\C_{\x^\jmath}$ and $\C_{\x^\kappa}$ increased while the error for $\C_{\x^\imath}$ decreased and therefore analysing the difference between the errors is an appropriate detection method for the power imbalance. This can be explained via \eqref{Error_power} and \eqref{error}. The decrease in the amplitude of $v_1$ and the increase in the amplitudes of $v_2$ and $v_3$ induced an increased $\imath$-imaginary part and decreased $\jmath$-imaginary and $\kappa$-imaginary parts of $x_l^{*}\left(n_1\right)x_l\left(n_2\right)$, so that the diagonalisation error of $\C_{\x^\imath}$ declined and the diagonalisation errors of $\C_{\x^\jmath}$ and $\C_{\x^\kappa}$ rose.
\section{Conclusion}\label{sect: conclusion}
Novel simultaneous matrix factorisation techniques for the covariance and complementary covariance matrices, which arise in widely linear quaternion algebra, have been introduced. The diagonalisation of quaternion symmetric matrices has been first addressed, as the Takagi factorisation of complex symmetric matrices cannot be readily extended to the quaternion domain. We have shown that the symmetry is insufficient for a quaternion matrix $\mathbf{A}$ to be diagonalisable, and instead the condition $(\mathbf{AA^*})^*=\mathbf{A}^*\mathbf{A}$ must also hold, which in general is not the case. The simultaneous diagonalisation of two quaternion covariance matrices, quaternion uncorrelating transform (QUT), has then been introduced as an analogue to the strong uncorrelating transform (SUT) in the complex domain. It has also been shown that for improper quaternion data, a typical case in quaternion signal processing applications, a single eigendecomposition of the covariance matrix is sufficient to diagonalise approximately the three complementary matrices. The analysis of the so introduced quaternion approximate uncorrelating transform (QAUT) has demonstrated its usefulness for improper quaternion signals, a typical case in practical applications. Simulation studies on synthetic and real-world signals support the proposed QUT and QAUT techniques.

\begin{appendices}
\renewcommand\thesection{\appendixname~\Alph{section}} 
\section{Properties of quaternion matrices}
For two general quaternion matrices, $\mathbf{A}$ and $\mathbf{B}$, the following properties hold \cite{CheongTook2011a}: 
\\ \newline
P1. $\mathbf{A}^{T*}=(\mathbf{A}^*)^T=(\mathbf{A}^T)^*~\text{and}~\mathbf{A}^{\alpha *}=(\mathbf{A}^*)^\alpha=(\mathbf{A}^\alpha)^*$, \\
P2. $\mathbf{A}^{\alpha T}=(\mathbf{A}^\alpha)^T=(\mathbf{A}^T)^\alpha ~ \text{and} ~ \mathbf{A}^{\alpha H}=(\mathbf{A}^\alpha)^H=(\mathbf{A}^H)^\alpha$, \\
P3. $(\mathbf{A}^\alpha)^\beta=(\mathbf{A}^\beta)^\alpha=\mathbf{A}^\gamma \quad \text{for distinct}~~\alpha,\beta,\gamma \in \{\imath,\jmath,\kappa\}$, \\
P4. $(\mathbf{AB})^*\neq\mathbf{A}^*\mathbf{B}^*$, \\
P5. $(\mathbf{AB})^T\neq\mathbf{B}^T\mathbf{A}^T$, \\
P6. $(\mathbf{AB})^H=\mathbf{B}^H\mathbf{A}^H$, \\
P7. $(\mathbf{AB})^\alpha=\mathbf{A}^\alpha\mathbf{B}^\alpha$, \\
P8. $(\mathbf{AB})^{\alpha H}=\mathbf{B}^{\alpha H}\mathbf{A}^{\alpha H}$.\\
\section{$\mathbb{C}^\alpha$-improper quaternion random vectors}
For distinct $\alpha,\beta \in \{\imath,\jmath,\kappa\}$, the Cayley-Dickson construction allows for the complex representation of a quaternion random vector $\mathbf{x}$, as
\vspace{-1mm}\[\mathbf{x}=\mathbf{z}_1+\mathbf{z}_2\beta
\]
where $\mathbf{z}_1$ and $\mathbf{z}_2$ are complex random vectors defined in the plane spanned by $\{1,\alpha\}$. The quaternion random vector $\x$ is $\mathbb{C}^\alpha$-improper if and only if both of the following two conditions are fulfilled \cite{Amblard2004}:
 
1) $\mathbf{z}_1$ and $\mathbf{z}_2$ are proper complex vectors, which is achieved when their real and imaginary parts are uncorrelated and with the same
variance,
 
2) $\mathbf{z}_1$ and $\mathbf{z}_2$ have different variances.
\section{Symmetric quaternion matrices}
Consider a $2 \times 2$ symmetric quaternion matrix $\mathbf{A}$, whereby \newline
\vspace{-2mm}\begin{align*}
\begin{split}
\mathbf{A}\mathbf{A}^*&=\begin{bmatrix}
a_{11} & a_{12} \\
a_{21} & a_{22} \\
\end{bmatrix} \begin{bmatrix}
a_{11}^* & a_{12}^* \\
a_{21}^* & a_{22}^* \\
\end{bmatrix} \\
&=\begin{bmatrix}
a_{11}a_{11}^*+a_{12}a_{21}^* & a_{11}a_{12}^*+a_{12}a_{22}^* \\
a_{21}a_{11}^*+a_{22}a_{21}^* & a_{21}a_{12}^*+a_{22}a_{22}^* \\
\end{bmatrix}
\end{split}
\end{align*}
\begin{align*}
\begin{split}
\mathbf{A}^*\mathbf{A}&=\begin{bmatrix}
a_{11}^* & a_{12}^* \\
a_{21}^* & a_{22}^* \\
\end{bmatrix} \begin{bmatrix}
a_{11} & a_{12} \\
a_{21} & a_{22} \\
\end{bmatrix} \\
&=\begin{bmatrix}
a_{11}^*a_{11}+a_{12}^*a_{21} & a_{11}^*a_{12}+a_{12}^*a_{22} \\
a_{21}^*a_{11}+a_{22}^*a_{21} & a_{21}^*a_{12}+a_{22}^*a_{22} \\
\end{bmatrix}
\end{split}
\end{align*}
By considering the symmetry $\mathbf{A}=\mathbf{A}^T$ ($a_{12}=a_{21}$), the matrices $(\mathbf{A}\mathbf{A}^*)^*$ and $\mathbf{A}^*\mathbf{A}$ can be structured as \\
\vspace{-2mm}\begin{align*}
\begin{split}
(\mathbf{A}\mathbf{A}^*)^* &=\begin{bmatrix}
|a_{11}|^2+|a_{12}|^2 & a_{12}a_{11}^*+a_{22}a_{12}^* \\
a_{11}a_{12}^*+a_{12}a_{22}^* & |a_{12}|^2+|a_{22}|^2 \\
\end{bmatrix} \\
\mathbf{A}^*\mathbf{A}~
&=\begin{bmatrix}
|a_{11}|^2+|a_{12}|^2 & a_{11}^*a_{12}+a_{12}^*a_{22} \\
a_{12}^*a_{11}+a_{22}^*a_{12} & |a_{12}|^2+|a_{22}|^2 \\
\end{bmatrix}
\end{split}
\end{align*}
\noindent Note that $(\mathbf{A}\mathbf{A}^*)^*=\mathbf{A}^*\mathbf{A}$ if their off-diagonal elements are commutative, implying that either the diagonal or off-diagonal elements of $\mathbf{A}$ are real-valued.
\end{appendices}
\section*{References} 
\bibliographystyle{elsarticle-num}
\bibliography{QuaternionJointDiagonalisation}
\newpage
\begin{figure}[h]
\centering
\subfloat
{\includegraphics[scale=0.3]{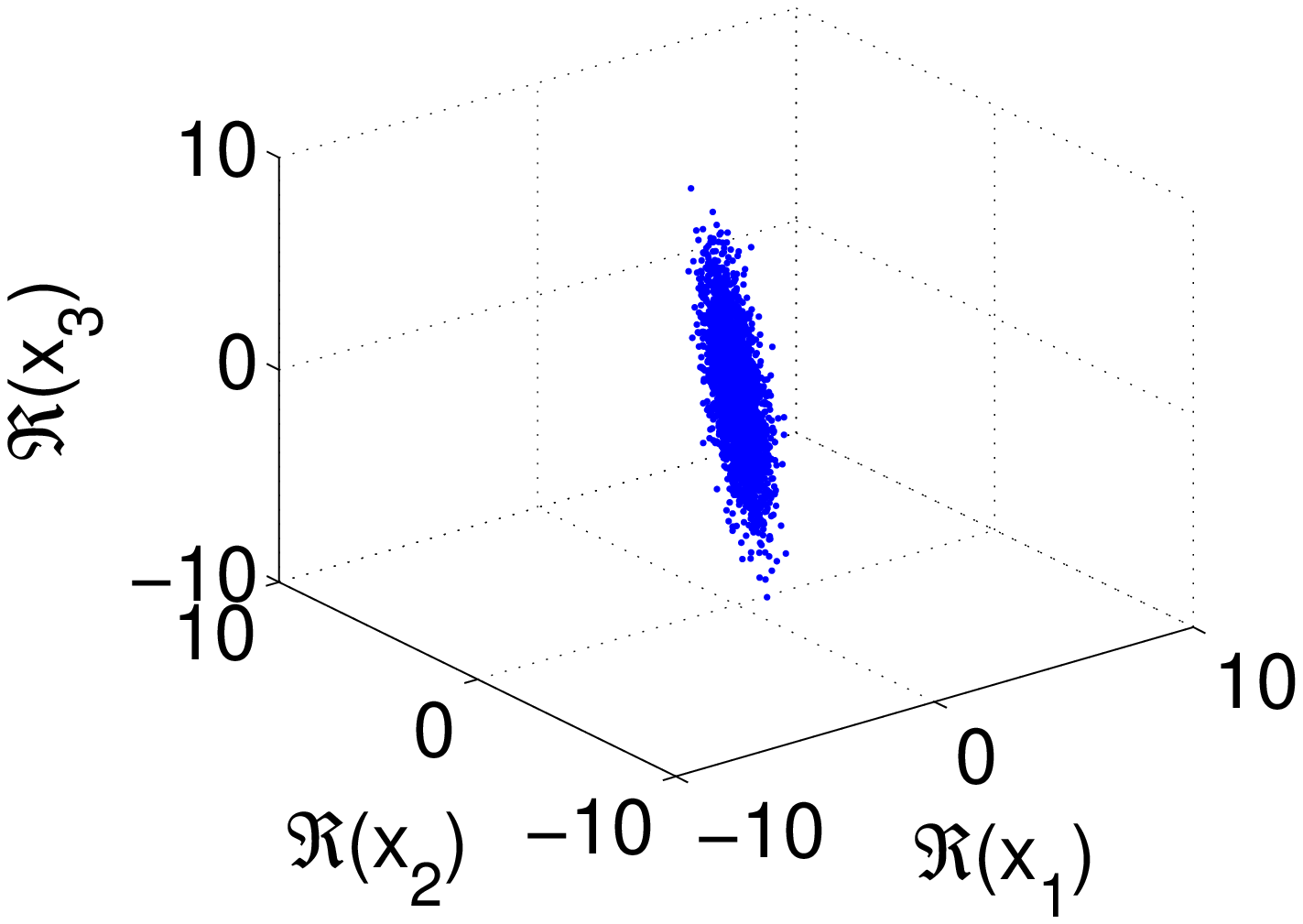}}
\subfloat
{\includegraphics[scale=0.3]{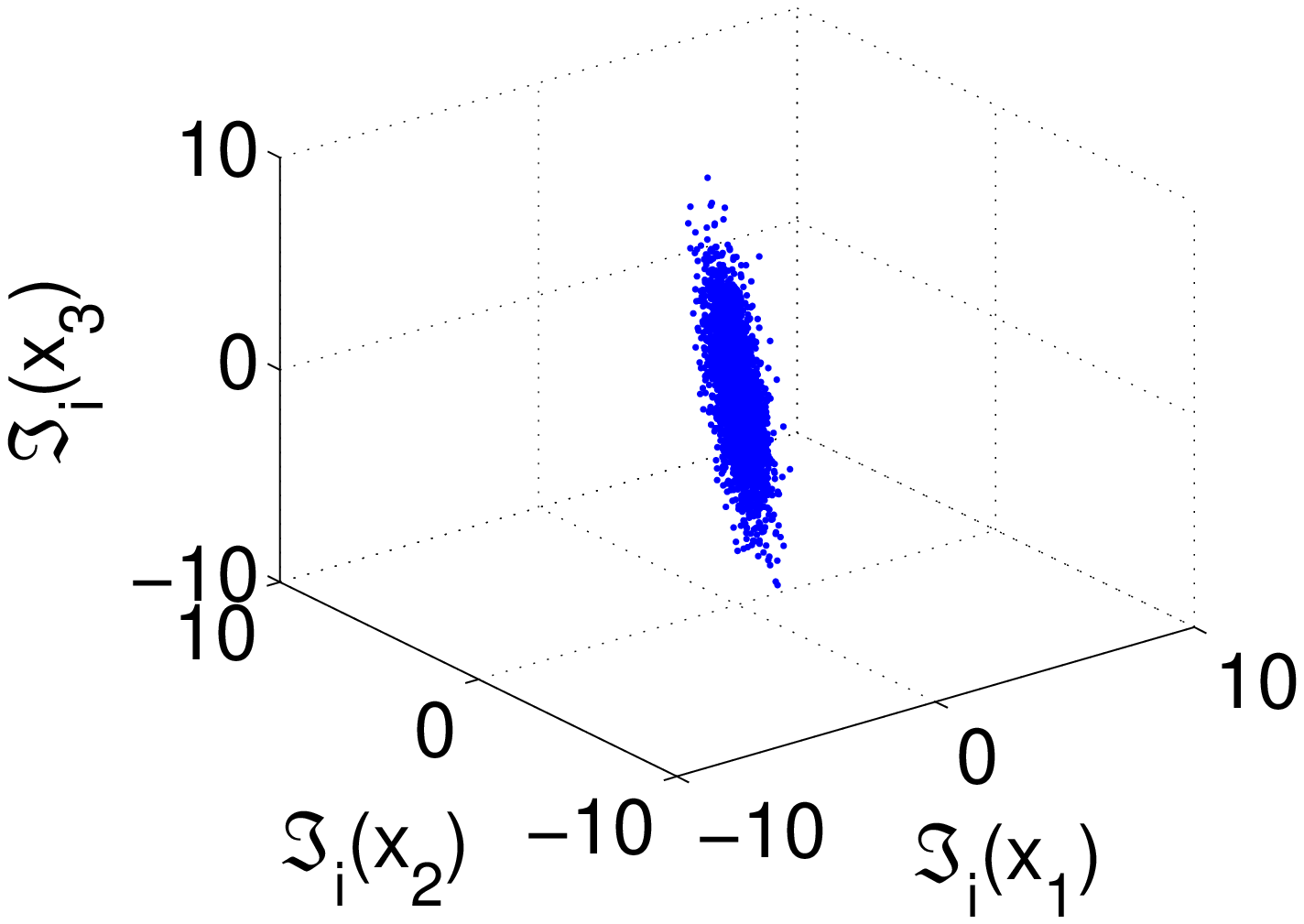}}

\subfloat
{\includegraphics[scale=0.3]{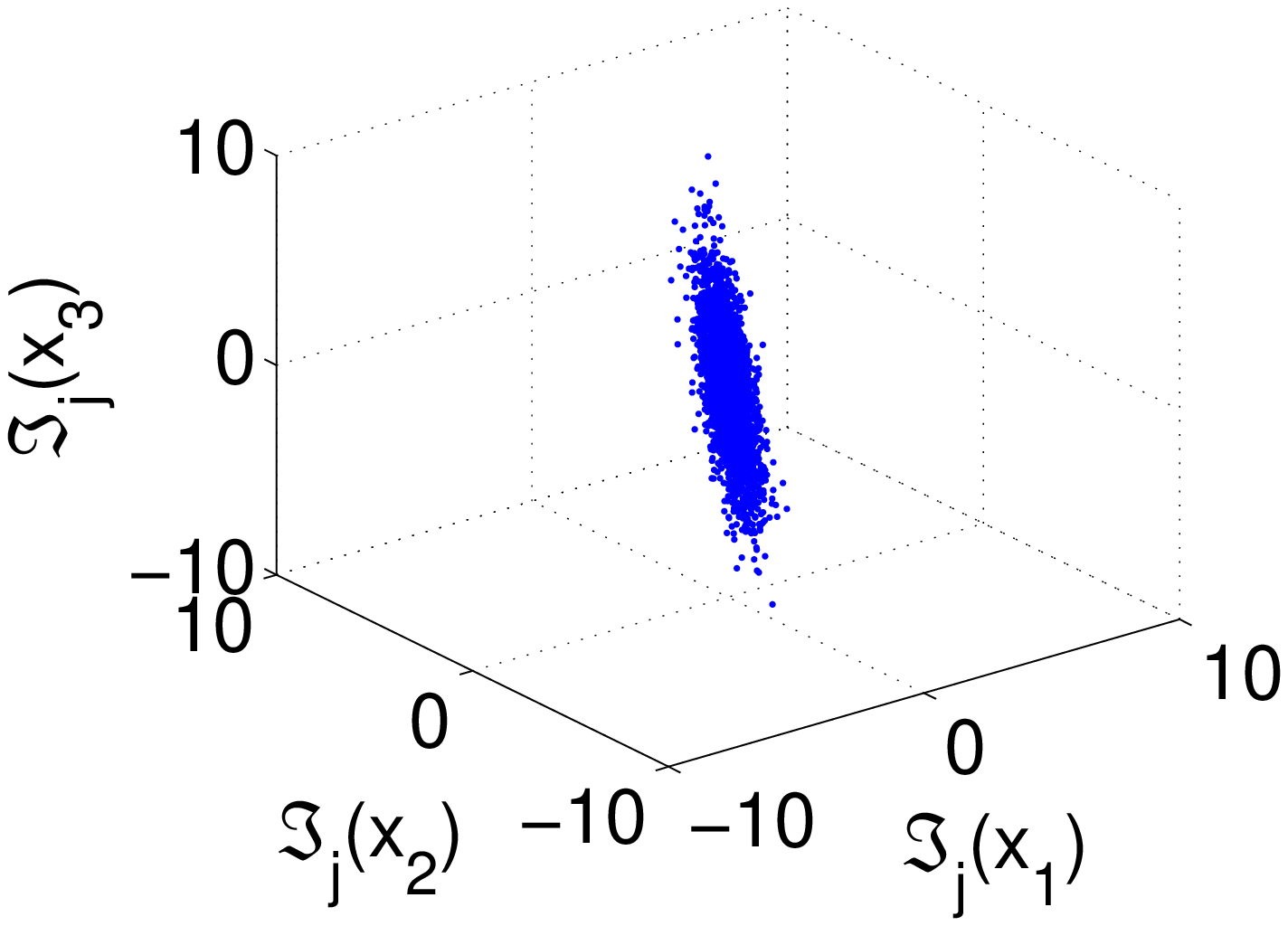}}
\subfloat
{\includegraphics[scale=0.3]{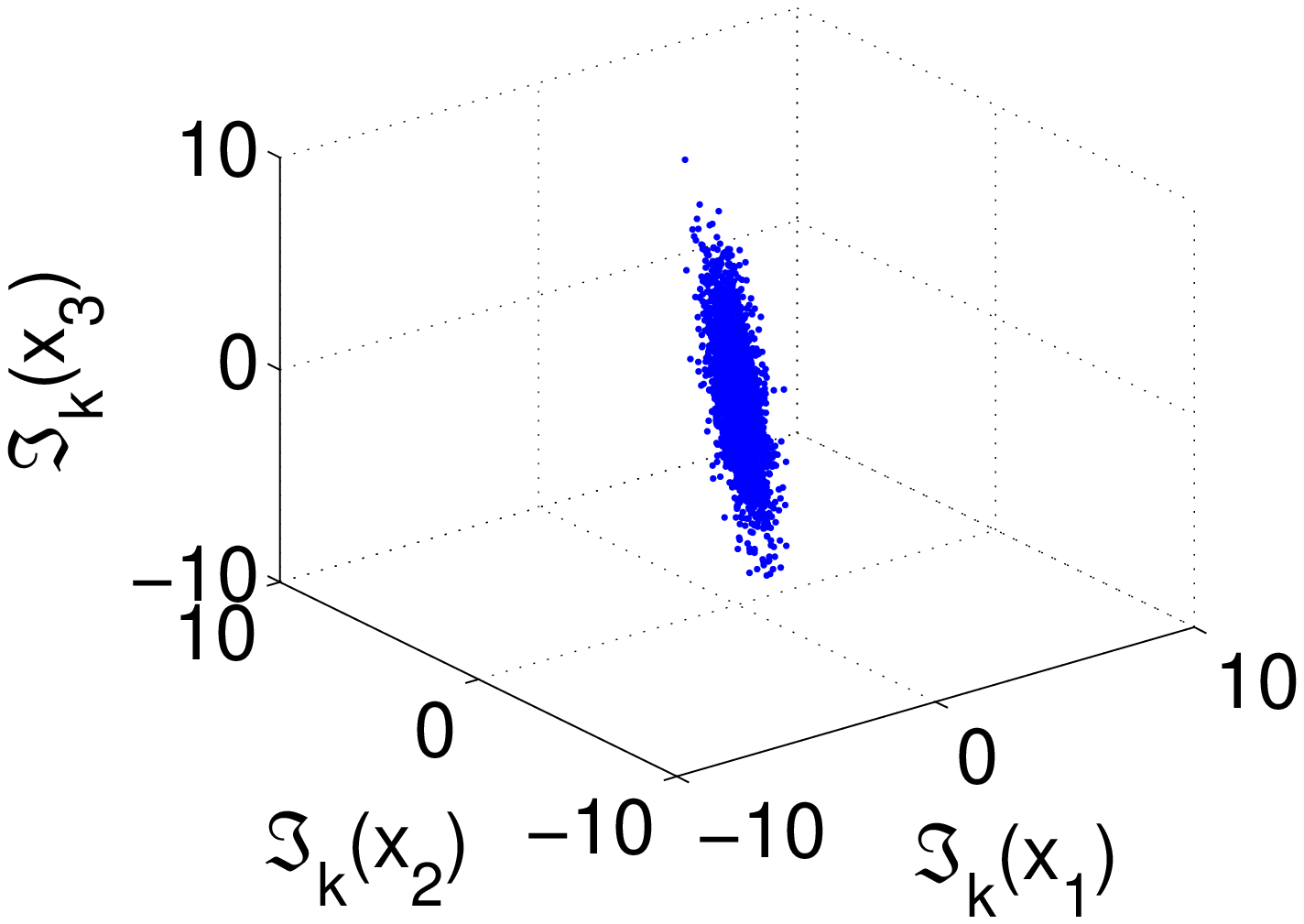}}
\caption{\label{fig:QUT1}Scatter plots, indicating the high correlation (improperness) of the three correlated $\mathbb{C}^\kappa$-improper quaternion signals, $\mathbf{x}_1$, $\mathbf{x}_2$ and $\mathbf{x}_3$.}
\end{figure}
\begin{figure}[h]
\centering
\subfloat
{\includegraphics[scale=0.3]{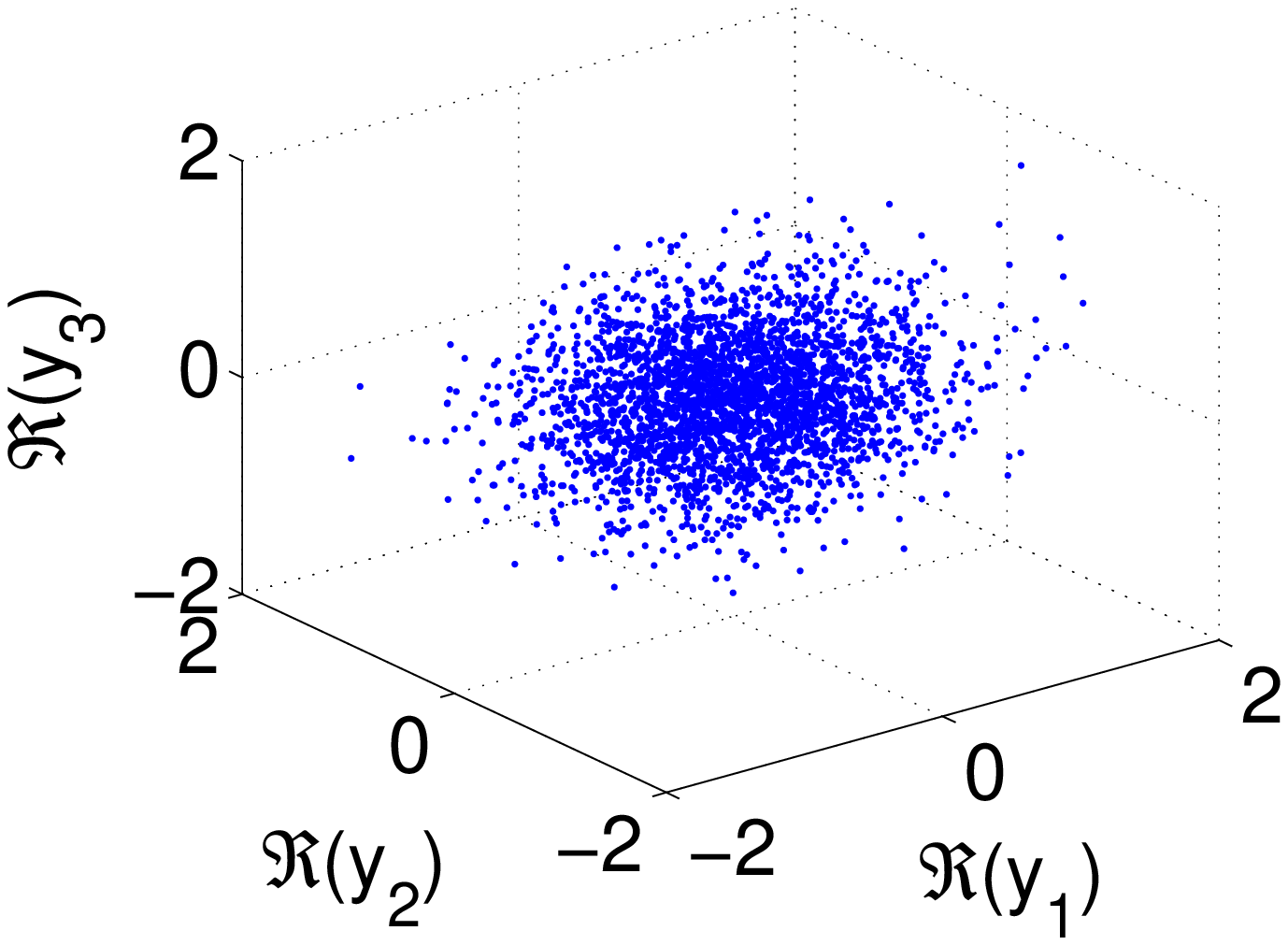}}
\subfloat
{\includegraphics[scale=0.3]{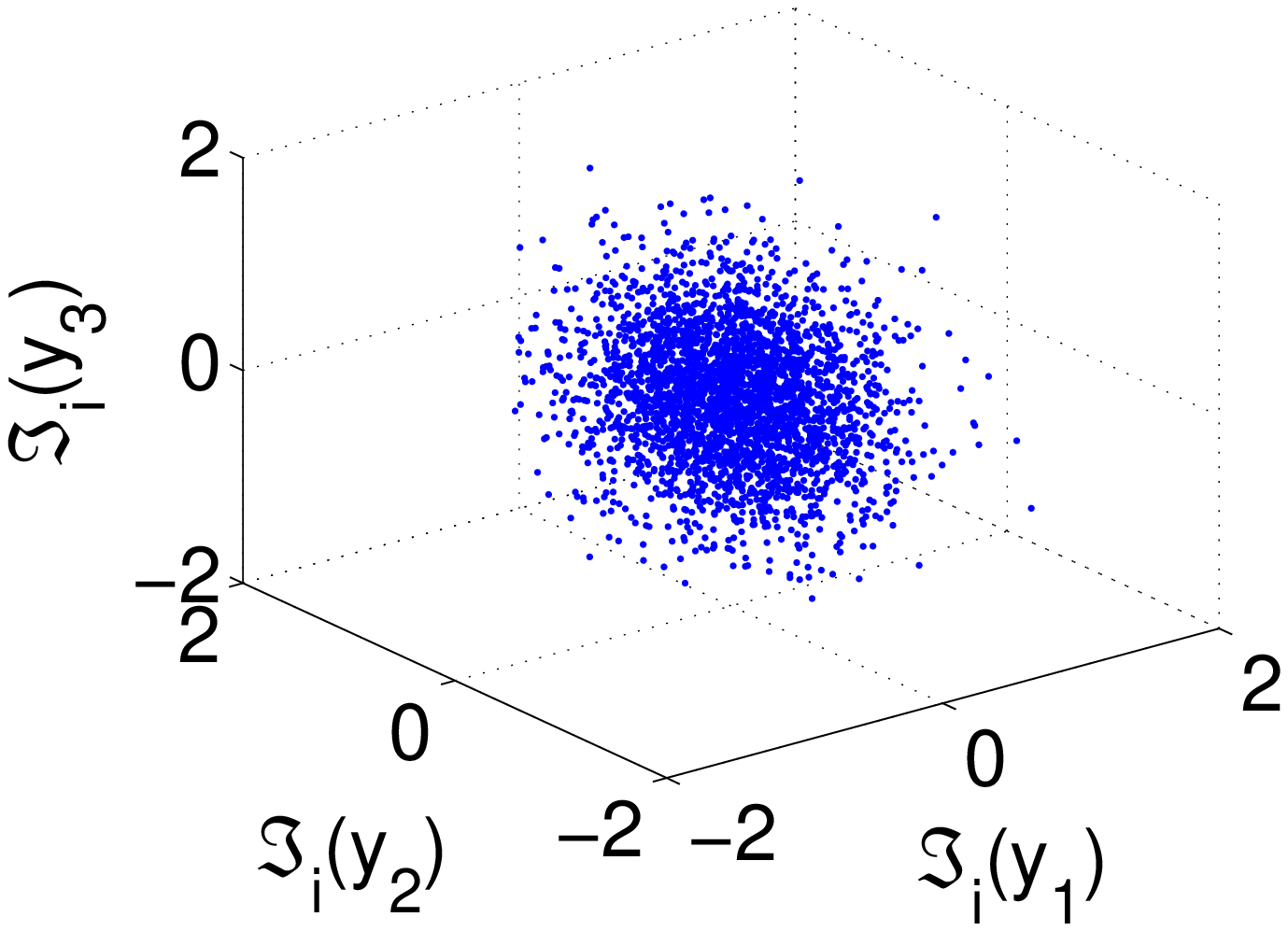}}

\subfloat
{\includegraphics[scale=0.3]{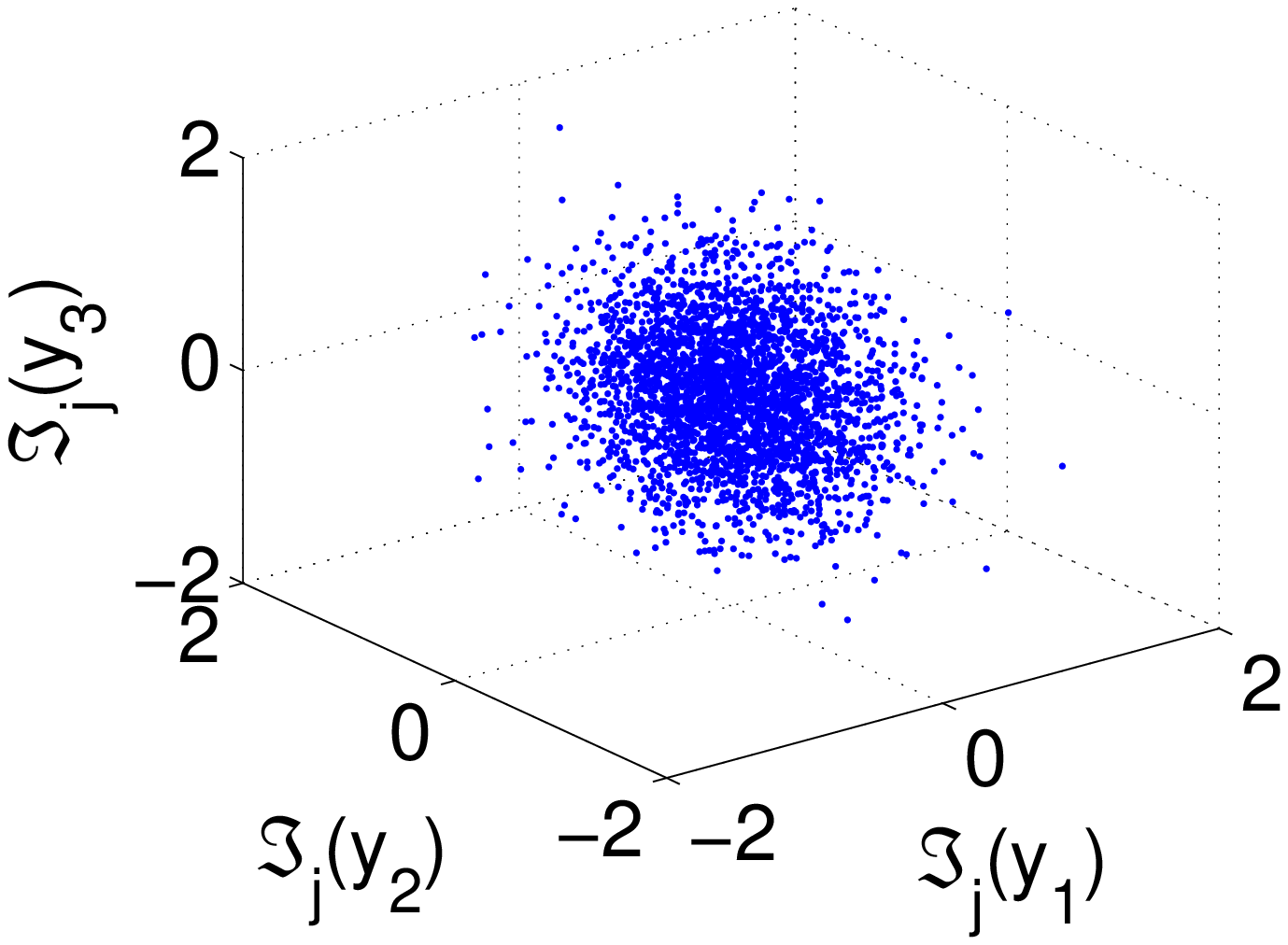}}
\subfloat
{\includegraphics[scale=0.3]{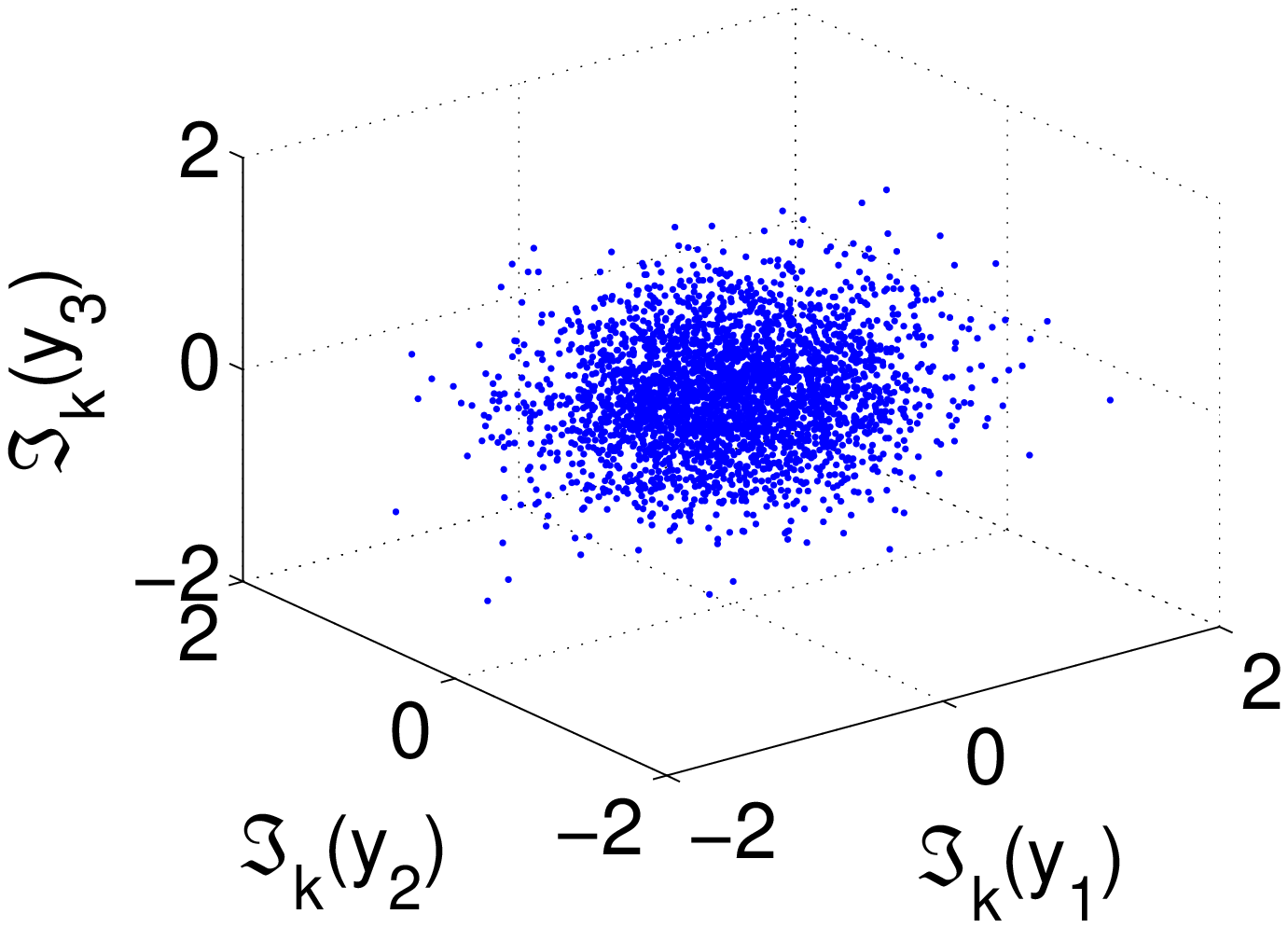}}
\caption{\label{fig:QUT2}Scatter plots, indicating the low correlation (properness) of the three $\mathbb{C}^\kappa$-improper quaternion signals, $\mathbf{y}_1$, $\mathbf{y}_2$ and $\mathbf{y}_3$, decorrelated via the QUT.}
\end{figure}

\begin{figure}[t]
\centering
\subfloat[Squared diagonal error of $\mathbf{C}_{\mathbf{x}^\imath}$]
{\includegraphics[scale=0.4]{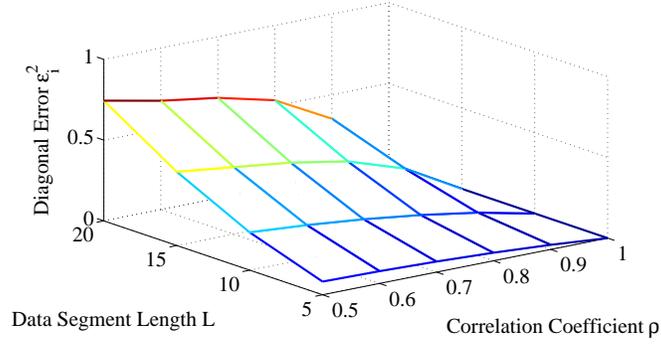}}
\hfil
\subfloat[Squared diagonal error of $\mathbf{C}_{\mathbf{x}^\jmath}$]
{\includegraphics[scale=0.4]{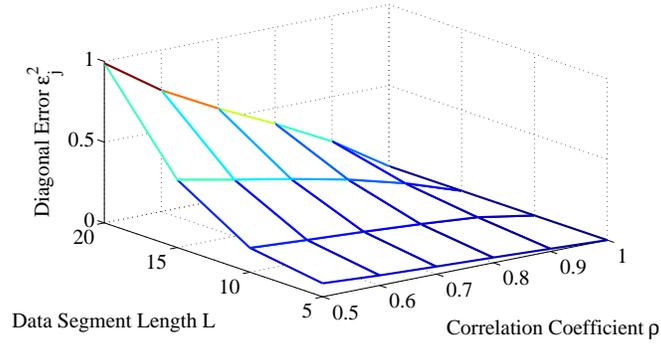}}
\hfil
\subfloat[Squared diagonal error of $\mathbf{C}_{\mathbf{x}^\kappa}$]
{\includegraphics[scale=0.4]{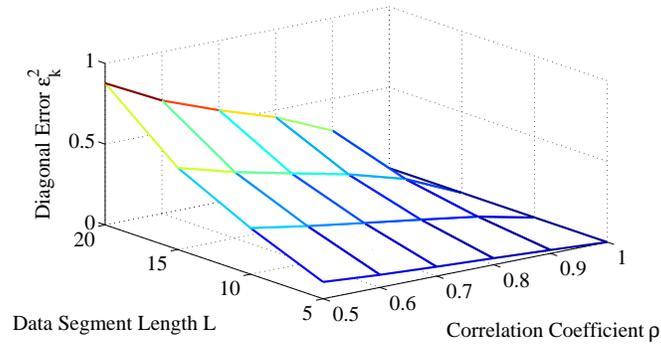}}
\caption{\label{fig:synthetic_univariate}Squared diagonal error [in $\%$] against the correlation degree and the data segment length for synthetic univariate data.}
\end{figure}

\begin{figure}[t]
\centering
\subfloat[Squared diagonal error of $\mathbf{C}_{\mathbf{x}^\imath}$]
{\includegraphics[scale=0.4]{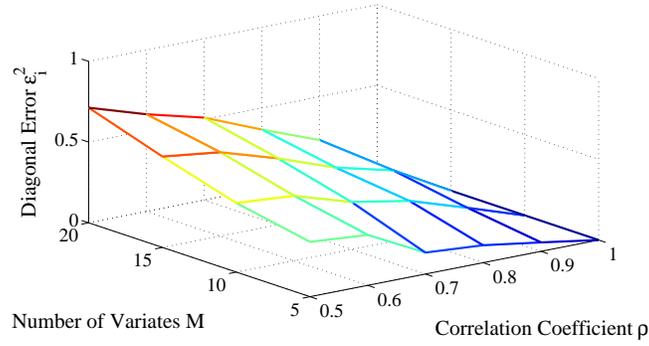}}
\hfil
\subfloat[Squared diagonal error of $\mathbf{C}_{\mathbf{x}^\jmath}$]
{\includegraphics[scale=0.4]{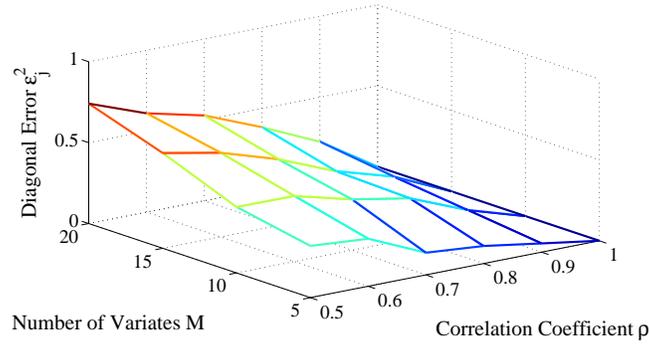}}
\hfil
\subfloat[Squared diagonal error of $\mathbf{C}_{\mathbf{x}^\kappa}$]
{\includegraphics[scale=0.4]{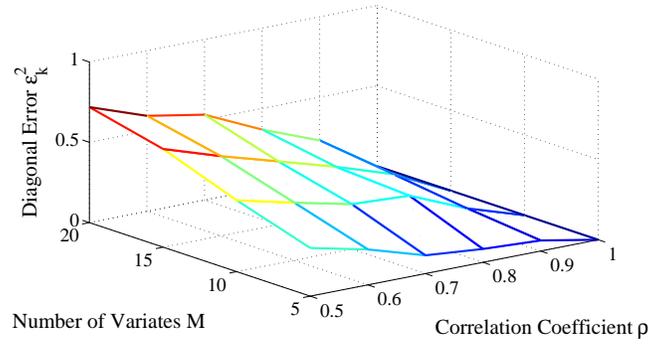}}
\caption{\label{fig:synthetic_multivariate}Squared diagonal error [in $\%$] against the correlation degree and the number of variates for synthetic multivariate data.}
\end{figure}

\begin{figure}[t]
\centering
\subfloat[Squared diagonal error of $\mathbf{C}_{\mathbf{x}^\imath}$]
{\includegraphics[scale=0.5]{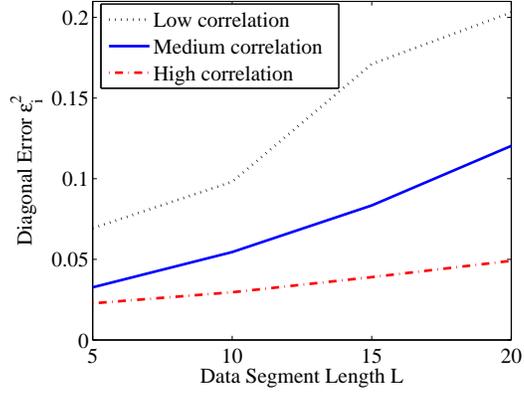}}
\vspace{-1mm}
\subfloat[Squared diagonal error of $\mathbf{C}_{\mathbf{x}^\jmath}$]
{\includegraphics[scale=0.5]{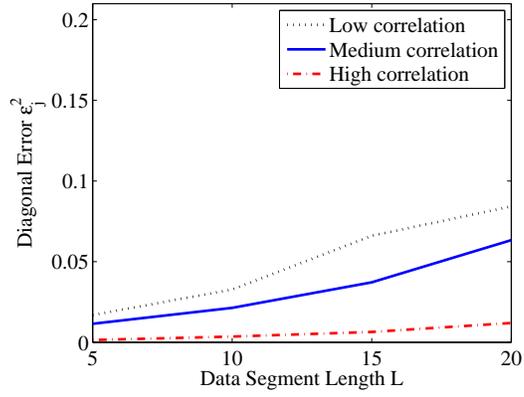}}
\vspace{-1mm}
\subfloat[Squared diagonal error of $\mathbf{C}_{\mathbf{x}^\kappa}$]
{\includegraphics[scale=0.5]{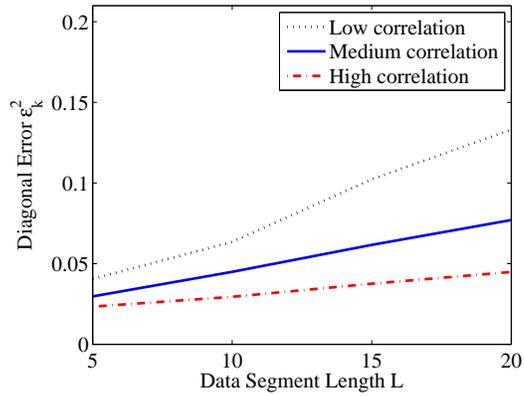}}
\vspace{-1mm}
\caption{\label{fig:EEG}Squared diagonal error [in $\%$] against the correlation degree and the data segment length for EEG data.}
\end{figure}

\begin{figure}[t]
\centering
{\includegraphics[scale=0.7]{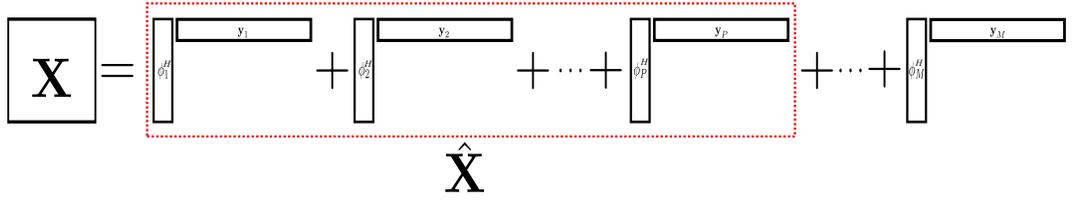}}
\caption{\label{fig:QAUTrankred}Rank reduction using the QAUT.}
\end{figure}

\begin{figure}[t]
\centering
{\includegraphics[scale=0.65]{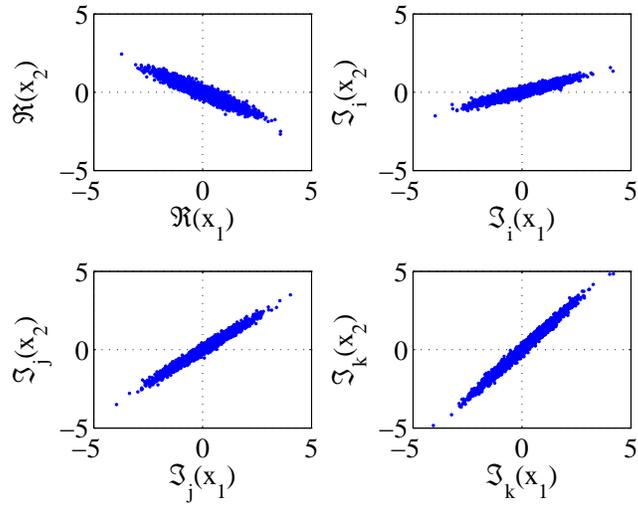}}
\caption{\label{fig:RankRed1}Scatter plots showing high correlation between the correlated $x_1$ and $x_2$ variables.}
\end{figure}

\begin{figure}[t]
\centering
{\includegraphics[scale=0.65]{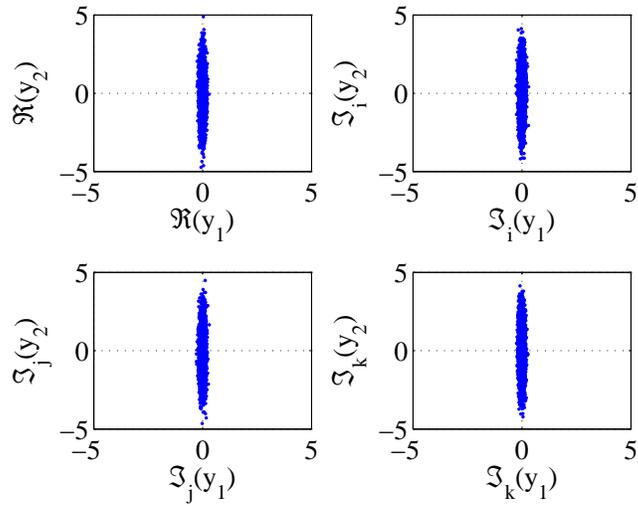}}
\caption{\label{fig:RankRed2}Scatter plots showing no correlation between the transformed $y_1$ and $y_2$ variables.}
\end{figure}

\begin{figure}[t]
\centering
{\includegraphics[scale=0.65]{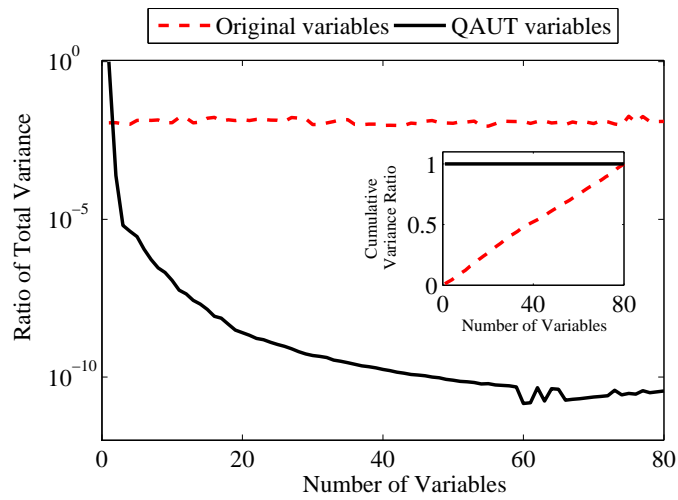}}
\caption{\label{fig:RankRedCorn}The logarithm of the ratio of the variance explained by each variable to the total variance. \emph{Insert diagram}: The cumulative ratio of variance against the number of variables for the original and decorrelated corn data.}
\end{figure}

\begin{figure}[t]
    \centering
    {\includegraphics[scale=0.52]{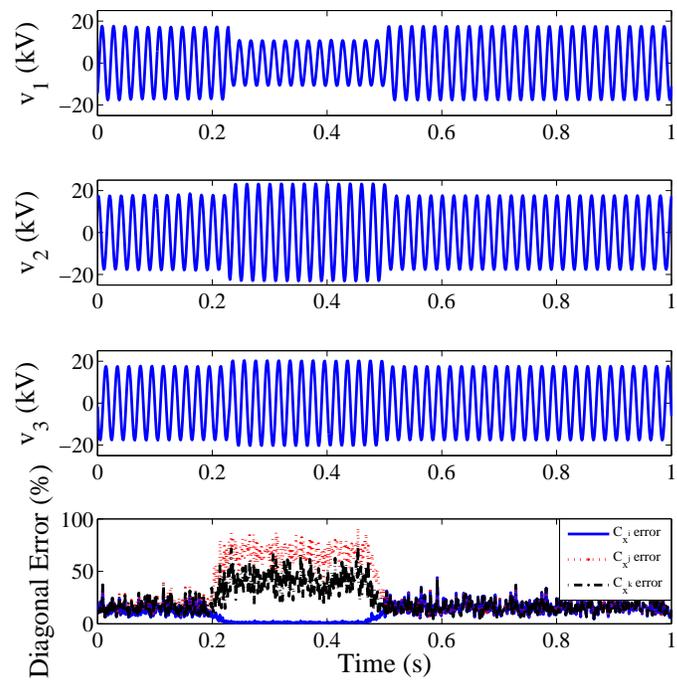}}
    \caption{Three-phase power signal and imbalance detection.}
    \label{fig:powimbal}
\end{figure}
\end{document}